\documentclass[a4paper, 10pt, reqno]{amsart}

\usepackage{ amssymb, amsmath, amsthm}
\usepackage{graphicx, psfrag, setspace, subfigure, gensymb}
\usepackage{amsfonts}
\usepackage{bbm}
\usepackage{fix-cm}
\usepackage{comment}
\usepackage{tikz,tikz-cd}
\usetikzlibrary{matrix,arrows,decorations.pathmorphing,decorations.pathreplacing}
\tikzset{commutative diagrams/diagrams={baseline=-2.5pt},commutative diagrams/arrow style=tikz}
\usepackage[colorlinks]{hyperref}
\usepackage{float}
\usepackage{setspace} 
\usepackage{enumitem}

\newcommand\Z{\mathbb Z}
\newcommand\C{\mathbb C}
\newcommand\R{\mathbb R}

\newcommand\A{\mathbb A}

\newcommand{\cC}{\mathcal{C}}
\newcommand{\cD}{\mathcal{D}}
\newcommand{\cE}{\mathcal{E}}
\newcommand{\cF}{\mathcal F}

\newcommand{\cO}{\mathcal{O}}

\newcommand{\cS}{\mathcal{S}}

\newcommand{\cW}{\mathcal{W}}

\newcommand\isoto{\stackrel{\sim}{\To}}
\newcommand\id{\mathrm 1}
\newcommand\acts{\curvearrowright}
\newcommand\into{\hookrightarrow}
\newcommand\To{\longrightarrow}

\newcommand\Hom{\operatorname{Hom}}

\renewcommand\P{\mathbb P}

\newcommand\Perf{\operatorname{Perf}}

\newcommand{\Crit}{\operatorname{Crit}}

\newcommand{\sslash}{/\!/}

\newcommand{\beq}[1]{\begin{equation}\label{#1} }
\newcommand{\eeq}{\end{equation}}
\newcommand{\pgap}{\vspace{5pt}}

\theoremstyle{plain}

\newtheorem{keyconj}{Conjecture}

\newtheorem{vagueconj}[keyconj]{``Conjecture''}

\theoremstyle{remark}
\newtheorem{rem}[equation]{Remark}

\newtheorem{remx}[equation]{Remark}
\newenvironment{remu}
  {\pushQED{\qed}\remx}
  {\popQED\endremx}

\theoremstyle{definition}

\newtheorem{eg}[equation]{Example}

\newtheorem{examplex}[equation]{Example}
\newenvironment{egu}
  {\pushQED{\qed}\examplex}
  {\popQED\endexamplex}

\makeatletter \@addtoreset{equation}{section} \makeatother

\let\oldtocsection=\tocsection
\let\oldtocsubsection=\tocsubsection
\let\oldtocsubsubsection=\tocsubsubsection
\renewcommand{\tocsection}[3]{\hspace{0em}\oldtocsection{#1}{#2}{#3}}
\renewcommand{\tocsubsection}[3]{ \hspace{1em} \oldtocsubsection{#1}{\small{#2}}{\small{#3}} }
\renewcommand{\tocsubsubsection}[3]{\hspace{2em}\oldtocsubsubsection{#1}{\small{#2}}{\small{#3}}}

\let\oldmarginpar\marginpar
\renewcommand\marginpar[1]{\-\oldmarginpar[\framebox{\setstretch{\marginparstretch}\begin{minipage}{\marginparwidth}{\raggedleft\scriptsize #1}\end{minipage}}]{\framebox{\setstretch{\marginparstretch}\begin{minipage}{\marginparwidth}{\raggedright\scriptsize #1}\end{minipage}}}}

\newcommand{\aand}{\quad\quad\mbox{and}\quad\quad}

\newcommand\Fuk{\mathcal{F}}
\newcommand\MF{\operatorname{MF}}
\newcommand\Ob{Ob}

\newcommand\Spec{\operatorname{Spec}}
\newcommand\mW{\breve{W}}
\newcommand\mX{\breve{X}}

\newcommand\mZ{\breve{Z}}
\newcommand\mD{\breve{D}}
\newcommand\mH{\breve{H}}
\newcommand\mw{\breve{w}}
\newcommand\mx{\breve{x}}
\newcommand\my{\breve{y}}
\newcommand\mz{\breve{z}}
\newcommand\um{\breve{u}}
\newcommand\mv{\breve{v}}
\renewcommand\mp{\breve{p}}
\newcommand\mq{\breve{q}}

\let\oldmarginpar\marginpar
\renewcommand\marginpar[1]{\-\oldmarginpar[\raggedleft\footnotesize #1]%
	        {\raggedright\footnotesize #1}}

\newskip\stdskip                      
\newcommand\mma{a}	
\newcommand\mmb{b}

\begin{document}

 \title{Equivariant Fukaya categories at singular values}
\author{Yank{\i} Lekili and Ed Segal}

\maketitle
\begin{abstract}

Given a Hamiltonian torus action on a symplectic manifold, Teleman and Fukaya have proposed that the Fukaya category of each symplectic quotient should be equivalent to an equivariant Fukaya category of the original manifold. We lay out new conjectures that extend this story - in certain situations - to singular values of the moment map. These include a proposal for how, in some cases, we can recover the non-equivariant Fukaya category of the original manifold starting from data on the quotient. 

To justify our conjectures we pass through the mirror and work out numerous examples, using well-established heuristics in toric mirror symmetry. We also discuss the algebraic and categorical structures that underlie our story.





\end{abstract}

\tableofcontents

\section{Introduction}\label{sec.intro}

We begin with an example that motivated our investigations. 

\begin{examplex} \label{eg.firsteg} Consider the exact symplectic manifold $T^*S^2$ with its canonical symplectic structure, which we can identify with the affine variety 
$$X = (xy+z^2=1)\subset \C^3$$
 with symplectic form restricted from $\C^3$. Projecting onto the $z$ direction exhibits $X$ as $\C^*$-fibration over $\C$, with two degenerate fibres over $\{\pm 1\}$. Write $S$ for the zero-section in $T^*S^2$, which is also the Lagrangian matching sphere over an arc $\Lambda\subset \C$ connecting the two critical values. 

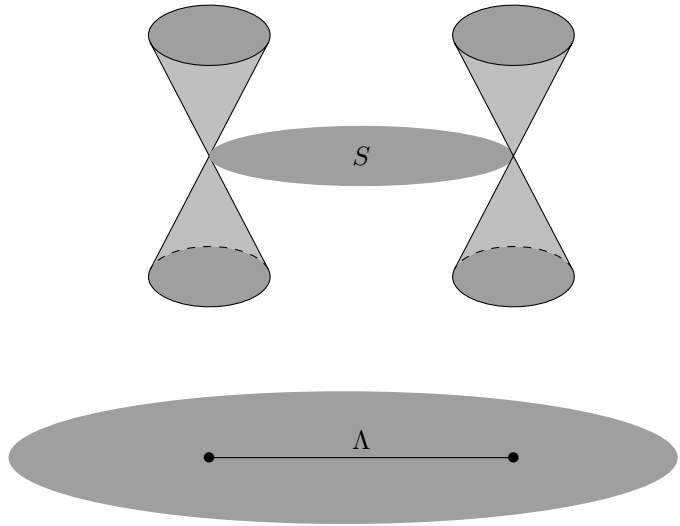
\begin{figure}[h!]
\begin{tikzpicture}[parabola/.style={very thin, red}, hyperbola/.style={very thin, violet}, circle/.style={very thin, blue}, scale=0.8]
\def\b{2}
\def\h{2}
\def\p{0.5}
\pgfmathsetmacro{\rx}{\b/2}
\pgfmathsetmacro{\ry}{\rx*\p}
\pgfmathsetmacro{\ta}{90-atan2(\h,\ry)}

\begin{scope}[xshift=2cm]
\fill[gray!50]
(0, \h) -- (\ta:\rx+0 and \ry) arc (\ta:180-\ta:\rx+0 and \ry) -- cycle;
\fill[gray!75] coordinate (bottom) ellipse [x radius=\rx, y radius=\ry];
\draw[dashed] (\ta:\rx+0 and \ry) arc (\ta:180-\ta:\rx+0 and \ry);
\draw (0, \h) -- (\ta:\rx+0 and \ry) arc (\ta:-180-\ta:\rx+0 and \ry) -- coordinate[pos=.4](hypbend) cycle;
\begin{scope}[rotate around={180:(0,\h)}]
\fill[gray!50]
(0, \h) -- (\ta:\rx+0 and \ry) arc (\ta:180-\ta:\rx+0 and \ry) -- cycle;
\fill[gray!75] coordinate (top) ellipse [x radius=\rx, y radius=\ry];
\draw (\ta:\rx+0 and \ry) arc (\ta:180-\ta:\rx+0 and \ry) (0, \h) -- coordinate[pos=.3](parabend) (\ta:\rx+0 and \ry) arc (\ta:-180-\ta:\rx+0 and \ry) -- cycle;
\end{scope}
\end{scope}

\fill[gray!75] (4.5,2) ellipse [x radius=\rx+1.5, y radius=\ry];

\begin{scope}[xshift=7cm]
\fill[gray!50]
(0, \h) -- (\ta:\rx+0 and \ry) arc (\ta:180-\ta:\rx+0 and \ry) -- cycle;
\fill[gray!75] coordinate (bottom) ellipse [x radius=\rx, y radius=\ry];
\draw[dashed] (\ta:\rx+0 and \ry) arc (\ta:180-\ta:\rx+0 and \ry);
\draw (0, \h) -- (\ta:\rx+0 and \ry) arc (\ta:-180-\ta:\rx+0 and \ry) -- coordinate[pos=.4](hypbend) cycle;
\begin{scope}[rotate around={180:(0,\h)}]
\fill[gray!50]
(0, \h) -- (\ta:\rx+0 and \ry) arc (\ta:180-\ta:\rx+0 and \ry) -- cycle;
\fill[gray!75] coordinate (top) ellipse [x radius=\rx, y radius=\ry];
\draw (\ta:\rx+0 and \ry) arc (\ta:180-\ta:\rx+0 and \ry) (0, \h) -- coordinate[pos=.3](parabend) (\ta:\rx+0 and \ry) arc (\ta:-180-\ta:\rx+0 and \ry) -- cycle;
\end{scope}
\end{scope}

\fill[gray!75] (4.2,-3) ellipse [x radius=\rx+4.5, y radius=\ry+0.6];

\node at (2,-3) {$\bullet$};
\node at (7,-3) {$\bullet$};

\node at (4.5,2) {$S$};

\draw (2,-3) -- (7,-3);
\node at (4.5,-2.7) {$\Lambda$};
\end{tikzpicture}
\caption{Cotangent bundle of two-sphere.}
\end{figure}

There is a Hamiltonian $S^1$ action on $T^*S^2$ given by rotating the fibres, \emph{i.e.}~rotating the phases of $x$ and $y$ in opposite directions. We are interested in the $S^1$-equivariant wrapped Fukaya category of $T^*S^2$ which we denote by  $\cW_{S^1}\!(T^*S^2)$. 

The construction of equivariant Fukaya categories with respect to an action of a Lie group $G$ is explained by Daemi and Fukaya in \cite[Theorem 8.3]{DF} \footnote{Other recent studies of $G$-equivariant Floer theory can be found in \cite{KLZ}, \cite{HLS}, \cite{Cazassus}. Though, the technical aspects of wrapped Floer theory involving non-compact invariant Lagrangians seems not to have been addressed yet.}. For this preliminary discussion we only need to know that an exact compact $S^1$-invariant Lagrangian is an object of this category, with endomorphisms given by $S^1$-equivariant Floer cohomology $HF^*_{S^1} (L,L)$,  which we can identify with the usual $S^1$-equivariant cohomology $H^*_{S^1}(L)$ (just as the Floer cohomology of an exact compact Lagrangian $HF^*(L,L)$ is isomorphic to $H^*(L)$).

The exact Lagrangian sphere $S\subset X$ is preserved by the $S^1$-action and, indeed, it defines an object of this category. The $S^1$-action on $S$ is the usual $S^1$-action on $S^2$ with two fixed points, so we can calculate that:
$$HF^*_{S^1}(S, S) = H^*_{S^1}(S) = \C[x,y]/xy$$
where $\deg x = \deg y =2$. Note that here we're using the (unique) grading structure on $T^*S^2$.

Now consider the set of regular values of our projection:
$$\mathcal{P} = \C\setminus \{\pm 1\}$$
This is the pair-of-pants. The wrapped Fukaya category $\cW(\mathcal{P})$ can easily be calculated, in particular the endomorphisms of the arc $\Lambda$ - which is a non-compact Lagrangian in $\mathcal{P}$ - is:
$$HW^*(\Lambda, \Lambda) = \C[x, y]/xy $$
Moreover, for an appropriate (and unique) grading structure on $\mathcal{P}$, we can arrange that $\deg x= \deg y=2$.

So we have observed that:
$$ HW^*(\Lambda, \Lambda)  = \C[x,y]/xy = HF^*_{S^1}(S, S) $$
This strongly suggests that there is quasi-equivalence of $\mathbb{Z}$-graded pre-triangulated categories:
\beq{eq.popeq} \mathcal{W}(\mathcal{P}) \simeq \mathcal{W}_{S^1}(T^*S^2) \eeq
We claim that this equivalence does indeed hold but we will not give a proof in this paper (although we will provide a heuristic mirror symmetry justification). Instead our focus will be on extrapolating from this example to some general conjectures. 

\end{examplex}

\subsection{Algebraic torus fibrations}

As a first generalization, consider the following set up. Let $Y= \mathbb{C}^n$, or a more general Liouville manifold,  and $f : Y \to \mathbb{C}$ be a holomorphic map with $0$ as a regular value. Then we can construct a fibration $\pi: X \to Y$ by setting $X$ to be the smooth hypersurface:
\[ X = \{ (x,y,\mathbf{w}) ;\, xy = f(\mathbf{w}) \} \; \subset \C^2\times Y \]
 The generic fiber of $\pi : X \to Y$ is isomorphic to $\C^*$ and it degenerates to a node along the smooth hypersurface $D = \{ f(\mathbf{w}) =0 \}$. The space $X$ admits an Hamiltonian $S^1$ action given by rotating the fibers: $ (x,y, \mathbf{w}) \mapsto (e^{i\theta}x,e^{-i\theta}y, \mathbf{w})$ for $e^{i\theta} \in S^1$.

More generally, suppose we have holomorphic line bundles $L_1,..., L_r$ on $Y$ and sections $f_i\in \Gamma(L_i^{\otimes 2})$. Assume that each hypersurface $D_i = \{ f_i=0 \}$ is smooth and that together they form a simple normal crossing divisor $D=\cup_{i=1}^r D_i$. Then we can form the smooth space:
\[ X = \{ (x_1,y_1,\ldots, x_r,y_r, \mathbf{w}) \,|\, x_iy_i = f_i(\mathbf{w}), \text{for all  } i \}\;\;\subset \operatorname{Tot}\left(\oplus_{i=1}^r L_i^{\oplus 2} \right) \]
The projection to $\mathbf{w}$ defines a fibration $\pi: X \to Y$ whose generic fiber is isomorphic to $(\mathbb{C}^*)^r$, with singular fibres appearing over the divisor $D$. We also have a Hamiltonian action of an $r$-dimensional torus $T$ on $X$ by rotating the fibers.

If $r>\dim Y$ then the intersection of all the $D_i$ is empty, and it follows that the diagonal $U(1)\subset T$ acts freely. Then the quotient $X' = X/\C^*$ is again a manifold fibering over $Y$, with generic fibre an algebraic torus, and singular fibres appearing over the same divisor $D$. It carries a Hamiltonian action of the rank $r-1$ torus $T/U(1)$. In fact if $r=\dim + k$ then there is a rank $k$ subtorus of $T$ which acts freely and we can quotient by any subtorus of it.
 
We'll refer to these $X\to Y$ or $X'\to Y$, with their Hamiltonian torus actions, as \emph{algebraic torus fibrations}.

\begin{keyconj}\label{conj.conics} 
Given an algebraic torus fibration $X\to Y$, we have a quasi-equivalence:
$$ \cW_{T}(X)_{-\mathbf{1}} \simeq  \cW(Y\setminus D)  $$
\end{keyconj}

Here if we work with $\mathbb{Z}$-graded categories, the grading on $Y \setminus D$ should be chosen so that it extends over $D$.  The meaning of the suffix $-\mathbf{1}$ will be explained in the next section; for the moment we ignore it.

If we delete the divisor $\pi^{-1}D$ from $X$ then what remains is just a principal $(\C^*)^r$ bundle over $Y\setminus D$. So it is not too surprising that there should be a quasi-equivalence:
$$\cW_T(X\setminus \pi^{-1}D) \simeq \cW(Y\setminus D)$$
and indeed this fits with a more general story about Hamiltonian reduction that we will discuss in the next section. However, $\cW_T(X)$ should be a deformation of $\cW_T(X\setminus \pi^{-1}D)$, since including the extra divisor will add terms to the $A_\infty$ structure. Conjecture \ref{conj.conics} is the claim that this deformation is in fact \emph{trivial}. 
\pgap

In Section \ref{sec.toricms} we will provide evidence for the conjecture using toric mirror symmetry; in some examples this amounts to an (excessively indirect!) proof. 

\begin{rem} There is also some evidence for our conjecture in the existing literature, not using mirror symmetry, when $X$ is a hypertoric variety associated to a special kind of hyperplane arrangement in $Y= \mathbb{C}^n$.

Fix a collection of points $p_1,..., p_k\in \C$. Taking the $n$th symmetric product produces a hyperplane arrangement  $Y \setminus D = \mathrm{Sym}^{n}(\mathbb{C} \setminus \{ p_1,\ldots, p_k\})$. In \cite{LPsym} the category $\cW(Y\setminus D)$ was described explicitly, by computing the endomorphism algebra of a particular generator. 

There is also a hypertoric variety $X$ which can be constructed as an algebraic torus fibration $\pi:X \to Y$ degenerating over $D$, as in Conjecture \ref{conj.conics}. The paper \cite{BLPW} studies certain algebras associated to $T$-invariant Lagrangians in $X$, which are expected to agree with their endomorphism algebras in $\cW_T(X)$. And, as was noted in \cite{LLM}, they are the same algebras as were found in \cite{LPsym}. 
\end{rem}

\subsection{Hamiltonian reduction}\label{sec.Hamiltonianreduction}

Suppose we have a Hamiltonian $S^1$ action on a symplectic manifold $X$. The information of the $S^1$ action appears in the Fukaya category as an invertible element:
$$s \in SH^*(X)$$
This gives rise to a natural automorphism of the identity functor, so for each object $L\in \cW(X)$ it provides an automorphism $s_{L}: L\isoto L$. Precisely, $s_{L} = \mathcal{CO}^0(s) \in HF^0(L,L)$ where $\mathcal{C}\mathcal{O}^0$ is the zeroth order part of the unital algebra map $\mathcal{CO}: SH^*(X) \to HH^*(CF^*(L,L))$ called the closed-open string map.    
 
This $s$ was originally constructed by Seidel  \cite{SeidelPi1} for compact $X$, and as an invertible element in the quantum cohomology $QH^*(X)$.  When $X$ is exact Seidel's construction yields an element in symplectic cohomology $SH^*(X)$ instead \cite{ritter}, but in either case we can map it to Hochschild cohomology using the closed-open map.\footnote{In the exact case the closed-open map $SH^*(X)\to HH^*(\cW(X))$ is known to be an isomorphism by work of Ganatra \cite{ganatra}.}



\pgap

In the previous section we discussed the equivariant wrapped category $\cW_{S^1}\!(X)$. In fact, for any given $\lambda\in\C^*$ one can construct a version of the equivariant category
$$\cW_{S^1}\!(X)_\lambda$$
starting from those objects of $\cW(X)$ such that $s_L= \lambda 1_{L}$. Teleman \cite{Teleman} refers to these categories as the `spectral components' of the equivariant Fukaya category.

 For example, an $S^1$-invariant Lagrangian $L$ (which is monotone, has minimal Maslov at least 2, and is equipped with a spin structure) provides an object of $\cW_{S^1}\!(X)_{\pm 1}$. But if we give $L$ a non-trivial local system, whose monodromy along $S^1$ orbits is $\lambda$, then we have an object of $\cW_{S^1}\!(X)_{\pm\lambda}$. 

In this paper we will focus almost exclusively on these spectral component categories, ignoring the global category $\cW_{S^1}(X)$. The reasons for this will become apparent in later sections. 

\begin{remu}\label{rem.signspin}
The sign ambiguity above comes from the choice of spin structure on $L$. Namely, along the $S^1$-orbit of a point on $L$, we have two trivialisation of the tangent bundle $TL$, one induced from the $S^1$ action and the other coming from the spin structure, and a negative sign appears if these do not agree, see \cite{Tonkonog1}. In particular modifying the local system by multiplying the monodromy by $-1$ is equivalent to changing the spin structure.

In Example \ref{eg.firsteg} our Lagrangian $S$ was a 2-sphere, hence it carries a unique spin structure which induces the bounding spin structure on each $S^1$ orbit. This makes it an object of $\cW_{S^1}\!(T^*S^2)_{-1}$. 

\end{remu}

Now consider the Hamiltonian reduction $ X\sslash_{\mma} \,S^1$ at some regular value $\mma\in \R$ of the moment map $\mu$. There is a Lagrangian correspondence
\begin{align} \label{momentLag} \Gamma = \{ (x, [x] ), \; \mu(x)=\mma \} \; \subset X^{-}\times  (X\sslash_{\mma} \,S^1) \end{align}
which induces a functor
\begin{align} \label{functor} \cW(X) \to \cW(X\sslash_{\mma} \,S^1) \end{align}
  via the theory of holomorphic quilts \cite{MWW, LL, Fuk}. Moreover, since $\Gamma$ is $S^1$-invariant we can use it to define a functor on the equivariant Fukaya category of $X$. Teleman conjectures that this gives an equivalence
\beq{eq.FukayaTeleman}  \cW_{S^1}\!(X)_{\pm e^{\mma}} \cong \cW(X\sslash_{\mma} \,S^1)\eeq
between the Fukaya category of the Hamiltonian reduction and the corresponding spectral component of the equivariant category.  More generally we can consider the spectral component at $\lambda = \pm e^{\mma+ i\mmb}$, and this will give the wrapped category of the quotient equipped with a B-field (or perhaps a more general bulk deformation). A theorem along these lines, for compact $X$, has been announced by Fukaya \cite{FukayaMiami}.

\begin{rem}
In general, one has to worry about whether $\Gamma$ is unobstructed. If $X$ is compact, it can occur that $\Gamma$ is only weakly unobstructed, which then necessitates choosing a bounding cochain on $\Gamma$ to make the functor (\ref{functor}) defined, in which case its target category may also get bulk deformed (as forced by the existence of the functor (\ref{functor}) \cite{Fuk}). 
\end{rem}

\begin{rem}\label{rem.exp(a)} The sign ambiguity in \eqref{eq.FukayaTeleman} is not present in Teleman's paper; we discovered it from examples but it can be explained as follows. Before we can use $\Gamma$ to define a functor we must first equip it with a spin structure, then it becomes an element of $\cW_{S^1}\!\big(X\times (X\sslash_\mma\, S^1)\big)_{\pm 1}$ with the sign determined as in Remark \ref{rem.signspin}. If $\Gamma$ is not simply connected then we may be able to produce two equivalences $\cW_{S^1}\!(X)_{e^\mma} \isoto \cW(X\sslash_\mathfrak
  {a}\,S^1)$ and $\cW_{S^1}\!(X)_{-e^\mma} \isoto \cW(X\sslash_\mma\,S^1)$ using different spin structures. In other examples there is no choice and one must just compute the sign. 
\end{rem}

Let us see what this point-of-view brings to our Example \ref{eg.firsteg}. 
The moment map there is $\mu = |x|^2 - |y|^2$. Any non-zero $\mma\in \R$ is a regular value of $\mu$, and produces the quotient $X\sslash_{\mma} \,S^1 \cong \C$. Since $\cW(\mathbb{C})\cong 0$ the corresponding spectral component should be zero. 

Given Remark \ref{rem.signspin}, our observation \eqref{eq.popeq} should really be the claim that:
$$\cW_{S^1}(X)_{-1} \cong \cW(P) $$
So the spectral component at $\lambda=-1$ is non-zero. This must correspond to the \emph{singular} value $\mma=0$ of the moment map, where we cannot do symplectic reduction. Instead, we are simply deleting the singularities of the moment map fibre $\mu^{-1}(0)$ and forming the quotient:
$$\mathcal{P} = \big(\mu^{-1}(0) - (0,0,\pm 1)\big) / S^1$$
With this prescription we are extending the equivalences
$$\cW_{S^1}(X)_{-e^\mma} \cong \cW(X\sslash_{\mma} \,S^1)$$
from \eqref{eq.FukayaTeleman} to the singular value $\mma=0$. 
 
From this example we draw the following general conclusion:

\begin{keyconj}\label{conj.singularfibres} Let $X$ be a Liouville manifold with a Hamiltonian $S^1$-manifold such that the fixed locus $X^{S^1}$ has codimension four, and assume there are no finite non-trivial stabiliser groups. Let $\mu: X \to \mathbb{R}$ be the moment map and let $\mma\in \R$ be a singular value in the interior of $\mathrm{Im}\mu$. Then there is a quasi-equivalence
$$\cW_{S^1}\!(X)_{- e^{\mma}} \cong \cW(U/S^1) $$
where $U$ is the smooth locus in $\mu^{-1}(\mma)$.
\end{keyconj}

As stated our conjecture may be too optimistic, as in general $\cW(U/S^1)$ may have to be bulk deformed.\footnote{When bulk deformations are involved, one usually deforms using a closed $2$-form (or a closed $l$-form in the ungraded case). It is known that $SH^*$ has an $L_\infty$-structure \cite{Siegel} and one can deform the wrapped Fukaya category along those $SH^*$ classes that satisfy the corresponding Maurer-Cartan equation. We also call these bulk deformations even though a geometric construction of these bulk deformations appears not to have been carried out in the literature. To deform the $A_\infty$-operations using a bulk class in $SH^2$, one should sum over disks with $r \geq 0$ interior punctures asymptotic to the bulk class, divided by $r!$.} However, we will see many examples below where such a bulk deformation is absent. 

\begin{rem} Following Remark \ref{rem.exp(a)} we should justify why we have written $-e^\mma$ instead of $e^\mma$.  The smooth locus $U$ is an $S^1$ bundle, but in $\mu^{-1}(\mma)$ some of the $S^1$ fibres collapse to points. So the Lagrangian correspondence $\Gamma$ should be given a spin structure which is the bounding spin structure on the orbits - hence the minus sign. This is essentially the same argument that we used for the $T^*S^2$ example in Remark \ref{rem.signspin}.
\end{rem}

\begin{rem}\label{rem.betanotzero}
The spectral components that are not covered either by Teleman--Fukaya's ideas or our Conjecture \ref{conj.singularfibres} are the components at $\lambda = - e^{\mma+ i\mmb}$ where $\mma$ is a singular value and $\mmb\neq 0$. The best we can say about these are that (i) they should be bulk deformations of $\cW(U/S^1)$, and (ii) they should be more similar to the components at generic values than they are to $\cW(U/S^1)$. From mirror symmetry we should expect to see special behaviour in complex codimension one. 
\end{rem}

Most of the discussion above generalizes easily to Hamiltonian actions of a higher rank torus $T=(S^1)^r$. Such an action produces $r$ Seidel elements $s_1,..., s_r$, and hence a category $\cW_T(X)_\lambda$ for each $\lambda \in (\C^*)^r$.  

For regular values of the moment map the corresponding component $\cW_T(X)_\lambda$ will be the Fukaya category of the corresponding symplectic quotient of $X$, possibly with a bulk deformation. At singular values, our conjecture is that we should instead - at least in some situations - take the quotient of the smooth part of the moment map fibre. 
\pgap

The set-up considered in Conjecture \ref{conj.conics} is a special case of this story.  A moment map value $(\mma_1,...,\mma_r)$ is regular if every $\mma_i$ is non-zero, and then the corresponding Hamiltonian reduction is the base $Y$. So taking the spectral component at the value $\lambda = (-e^{\mma_1},..., -e^{\mma_r})$ will produce $\cW_{S^1}\!(X)_\lambda \cong \cW(Y)$. At a more general $\lambda$, provided each $\lambda_i$ lies off the unit circle, we will get $\cW(Y)$, possibly with a bulk deformation.

But suppose $\lambda_1=-1$. Then the fibre of the moment map has singularities over the divisor $D_1$, and our conjecture is that that spectral component is 
$$\cW_{S^1}\!(X)_\lambda \cong \cW(Y\setminus D_1)$$
 bulk deformed by $(\mmb_2,..., \mmb_r)$.  If several $\lambda_i$ are equal to $-1$ then we must delete several divisors, and at the most special value 
$$ -\mathbf{1} = (-1,..., -1)$$
 we get the statement of Conjecture \ref{conj.conics}. 

\begin{egu}\label{eg.simplest0} The simplest possible example of Conjectures \ref{conj.conics} and \ref{conj.singularfibres} is to set $X=\C^2$ with the $S^1$ action $\theta: (x,y) \to (e^{i\theta}x, e^{-i\theta}y)$.  We can view it as an algebraic torus bundle $\pi:X \to Y=\C$ with $\pi(x,y)=xy$, having a single singular fibre over $D=\{0\}$. The Hamiltonian reduction at generic values of $\mu$ is $\C$,  and at the singular value $\mma=0$ the quotient of the smooth locus is $\C^*$. So we expect that
$$\cW_{S^1}(X)_{-e^\mma}\cong \cW(\C)\cong 0 $$
 if $\mma\neq 0$, but that $\cW_{S^1}\!(X)_{-1}\cong \cW(\C^*)$ which is not zero. 

Consider the invariant Lagrangian torus $L=\{|x|=|y|=1\}$, and equip it with the canonical spin structure and some choice of $\C^*$ local system. Then we have an object of $\cW_{S^1}\!(X)_\lambda$ where $\lambda$ is the monodromy along $S^1$ orbits in $L$. Following the discussion after Conjecture \ref{conj.conics} we could delete $\pi^{-1}D$ and consider just the $\C^*$ bundle $\pi: (\C^*)^2\to \C^*$. In this space  $L$ is exact, and since the $S^1$ action is free we can compute
$$HF_{S^1}(L,L) = HF(L/S^1, L/S^1) = H^*(S^1)$$
regardless of the local system. This is consistent with with the equivalence \eqref{eq.FukayaTeleman} that says
$$\cW_{S^1}\!(X\setminus \pi^{-1}D)_{-e^\mma} \cong \cW(\C^*)$$
for all values of $\mma$.  But if we do the computation in $X$ instead this answer gets deformed by holomorphic discs. Indeed the moduli space of discs bounded by $L$ is diffeomorphic to $S^1 \sqcup S^1$,  and the $S^1$ action is free on each component \cite[Lem.~2.16]{SeidelLES}. Moreover the boundary classes of the two $S^1$-families of discs differ by an orbit of a point in $L$,\footnote{See \cite[Lem.~2.19]{LM} for a generalisation of this computation.} and their contributions cancel each other exactly when $\lambda = -1$. 

So if $\lambda\neq -1$ then this object is isomorphic to zero in $\cW_{S^1}\!(X)_\lambda$, but if $\lambda=-1$ we have a whole family of non-zero objects parameterised by the remaining monodromy $\nu\in \C^*$ of the local system. This supports the claim that $\cW_{S^1}\!(X)_{-1}\cong \cW(\C^*)$ but that $\cW_{S^1}\!(X)_\lambda \cong 0$ for all other values of $\lambda$. It seems that the deformation 
$$\cW_{S^1}\!(X\setminus \pi^{-1}D)_\lambda \leadsto \cW_{S^1}\!(X)_\lambda$$
 is trivial exactly for the special value $\lambda =-1$.

Following Remark \ref{rem.betanotzero} we can look in particular at the components $\cW_{S^1}\!(X)_{-e^{i\mmb}}$ for $\mmb\neq 0$. These all appear to be zero, and hence non-trivial bulk deformations of $\cW_{S^1}\!(X)_{-1}$.

  In Example \ref{eg.simplest2} we will show that all these claims are consistent with mirror symmetry.

\end{egu}

In fact we haven't found any examples for Conjecture B, as stated, that go beyond the set-up of Conjecture A. But we do have a couple of rather speculative examples where we drop the condition that $X$ is Liouville - see Section \ref{sec.beyondconics}.

\subsection{Mirror symmetry} \label{sec.MSintro}

Given a Hamiltonian action $S^1\acts X$, the Seidel element makes the wrapped Fukaya category $\cW(X)$ - at least at the level of homology - linear over the ring $\C[s, s^{-1}]$. This can be understood as part of a more general story of `topological group actions' on categories envisaged by Teleman \cite{Teleman}, see Section \ref{sec.topologicalactions}.

Now suppose that $X$ is mirror to an algebraic variety $\mX$. Then, since $\cW(X)=D^b(\mX)$, the mirror to $s$ must be an invertible element $\sigma$ in:
\beq{HH0} HH^0(D^b(\mX)) = \Gamma(\cO_{\mX}) \eeq
If we have an $S^1$ action on the symplectic side then on the mirror we have a function $\sigma:\mX\to \C^*$. 

\begin{eg}\label{eg.seideloncylinder} Take the cylinder $X=T^*S^1$ with the obvious $S^1$ action. The mirror is $\mX=\C^*$, which we can view as the space of $\C^*$-local systems on the zero section $S^1\subset X$. If we equip this $S^1$ with the non-bounding spin structure then it follows that the mirror to the $S^1$ action is the identity function $\sigma: \mX \to \C^*$.
\end{eg}

More generally $X$ might be mirror to a Landau-Ginzburg model $(\mX, \mW)$.  Then $\cW(X)$ is equivalent to the category of matrix factorizations $\MF(\mX, \mW)$, but still a function $\sigma:\mX\to \C^*$ does provide a natural automorphism of this category, so a possible mirror to the $S^1$ action on $X$. 

\begin{vagueconj}\label{conj.mirror}
Suppose we have a Hamiltonian $S^1$ action on a Liouville manifold $X$. Suppose $X$ has a mirror Landau-Ginzburg model $(\mX, \mW)$, and that the $S^1$ action is mirror to a function $\sigma: \mX\to \C^*$. Then for every $\lambda\in \C^*$ we have an equivalence
$$\cW_{S^1}\!(X)_\lambda \;\cong\;  \MF\!\big(\mZ_\lambda,  \mW|_{\mZ_\lambda} \big) $$
where $\mZ_\lambda\subset \mX$ denotes the hypersurface $\sigma^{-1}(\lambda)$. 
\end{vagueconj}

The generalization to higher-rank torus actions is obvious.

This claim will be central to all the mirror symmetry evidence that we provide for our other conjectures. However, it is not really a precise conjecture because we have not specified what we mean by `mirror'. In particular it's not enough to just assume that $\cW(X)\cong \MF(\mX, \mW)$, as we shall see in Section \ref{sec.topologicalactions}. 

\begin{egu}\label{eg.nodedeforming1} Consider
$$X = \mathbb{C}^2 \setminus \{ xy=1 \}$$
equipped with the restriction of the standard symplectic form on $\C^2$. This is a log-Calabi Yau surface which is known to be self-mirror (\cite{pascaleff}). It is an important example because it is the simplest example of a symplectic manifold admitting a Lagrangian torus fibration that has a singular fibre. We write
$$ \mX = \mathbb{C}^2 \setminus \{ \mx\my=1 \} $$
for the mirror. 

 Now add a Hamiltonian $S^1$ action on $X$ by  $e^{i\theta}(x,y) = (e^{i\theta}x, e^{-i\theta}y)$.  Using some toric mirror symmetry (see Example \ref{eg.nodedeforming2}) one can argue that on $\mX$ this becomes the non-vanishing function:
  \[ \sigma =  -1-\mx\my \]
If $\lambda\in \C^*$ with $\lambda\neq -1$ then $\sigma^{-1}(\lambda)\cong\C^*$. So the claim is that $\cW_{S^1}\!(X)_\lambda\cong D^b(\C^*)$. But for $\lambda=-1$ we are claiming that $\cW_{S^1}\!(X)_{-1}$ should be equivalent to the derived category of the node $\mZ_{-1}= \{\mx\my=0\}$. 

We can also relate these spectral components to symplectic reduction as we discussed in the previous section. For a regular moment map value $\mma\in \R \setminus 0$,  Fukaya-Teleman tell us \eqref{eq.FukayaTeleman}  that 
$$\cW_{S^1}\!(X)_{-e^\mma}\cong \cW(X\sslash_{\mma}\, S^1) = \cW(\C \setminus 1)$$
which is consistent with our mirror symmetry claims since $\C\setminus 1$ is indeed mirror to $\C^*$.  At the singular value $\mma=0$ we can apply our Conjectures \ref{conj.conics} or \ref{conj.singularfibres} which tell us to delete the singularity from $\mu^{-1}(0)$ before we take the quotient. The result is $\mathbb{C} \setminus \{0, 1\}$ which is the pair-of-pants,\footnote{The grading on this pair-of-pants is different from the one in Example \ref{eg.firsteg}.} and this is indeed mirror to the node $\mZ_{-1}$.

We can also move off the real line and set $\lambda = -e^{\mma+i\mmb}$. If $\mma\neq 0$ this has no effect; all bulk deformations of $\C\setminus \{1\}$ are trivial (there are no $SH^2$ classes to deform with) and $\cW(\C\setminus \{ 1 \}) \cong D^b(\C^*)$, consistent with the mirror. But this is not true for $\mma=0$ and we can see from the mirror that $\cW_{S^1}\!(X)_{-e^{i\mmb}}$ should be a non-trivial deformation of $\cW_{S^1}\!(X)_{-1}$ (\emph{c.f.}~Remark \ref{rem.betanotzero} again).

\end{egu}

\begin{remu}\label{rem.notfibres} It is tempting to understand ``Conjecture'' \ref{conj.mirror} in the following way.
  Both $\cW(X)$ and $\MF(\mX, \mW)$ are categories linear over $\C[s, s^{-1}]$, and they are equivalent over this base. For any fixed value $s=\lambda\in \C$ we should be able to take the fibres on both sides, and get equivalent categories. 

This interpretation is a useful first approximation, but it suffers from the following two defects.
\begin{enumerate}\item The Seidel element is only central in $\cW(X)$ at the homology level, so one should be careful with the definition of `linear over $\C[s, s^{-1}]$'.  To be able to base-change $\cW(X)$ to individual values $s=\lambda$ presumably would require some $E_2$ structure. There is recent progress on constructing such $E_2$ structure in \cite{OT}. 

\end{enumerate}
More seriously:
\begin{enumerate}\addtocounter{enumi}{1}
\item It is \emph{not true} in general that the fibre of $\MF(\mX, \mW)$ at $s=\lambda$ is equal to the category of matrix factorizations on the hypersurface $\mZ_\lambda$. See Section \ref{sec.topologicalactions}.
\end{enumerate}
\end{remu}

Note that the two functions $\mW$ and $\sigma$ on the mirror $\mX$ are playing very different roles, and this becomes clearer when we discuss grading. To give $\cW(X)$ a $\Z$-grading we must choose a grading structure, which is a homotopy class of sections of $LGr(TX)$. On the mirror this usually\footnote{But not always, see \cite{SegalTwisted}.} corresponds to a choice of `R-charge', \emph{i.e.}~a $\C^*$ action on $\mX$ for which $\mW$ has weight 2. The function $\sigma$ on the other hand must have weight zero. 

The presence of R-charge alters the grading on Hochschild cohomology so destroys the equality \eqref{HH0}, even in the case $\mW=0$. So in principle the mirror to $S^1\acts X$ might be a $\sigma$ with some higher-order terms given by polyvector fields. But we will not encounter any such examples in this paper.

\begin{remu}\label{rem.equivariantmirror} We have only made claims about the mirrors to spectral components $\cW_{S^1}\!(X)_\lambda$, not about the equivariant category $\cW_{S^1}\!(X)$ itself. However, in the $\Z$-graded situation we believe that the spectral components, together with the non-equivariant category $\cW(X)$, should determine $\cW_{S^1}\!(X)$. This too will be explained in Section \ref{sec.topologicalactions}. So in the $\Z$-graded world, we can perhaps read ``Conjecture" \ref{conj.mirror} as a \emph{definition} of an `$S^1$-equivariant mirror'.\footnote{See the recent paper of Aganagic \cite{Aganagic} for some related ideas.}

 It would be very interesting to understand the mirror to $\cW_{S^1}\!(X)$ directly, and we hope to address this in future work.
\end{remu}

\subsection{Recovering the non-equivariant category}\label{sec.recovering}

Given $S^1\acts X$ we have for each $\lambda\in \C^*$ a spectral component $\cW_{S^1}\!(X)_\lambda$ of the equivariant Fukaya category. These categories have some important extra structure, they are linear over the ground ring of $S^1$-equivariant cohomology:
$$H^\bullet_{S^1}(pt) = \C[t], \quad \deg t =2 $$
This structure is built into the construction and all $A_\infty$ structure maps respect it. It is therefore possible to take the fibre of $\cW_{S^1}\!(X)_\lambda$ at $t=0$. The result is a subcategory 
$$\cW_{S^1}\!(X)_\lambda|_{t=0} \subset \cW(X) $$
of the ordinary wrapped Fukaya category of $X$, it is the full subcategory of objects $L$ with $s_{L}=\lambda 1_{L}$.
 
\begin{examplex} In Example \ref{eg.firsteg}, the endomorphism algebra of $S$ in $\cW_{S^1}\!(X)_{-1}$ is the equivariant cohomology $H^\bullet_{S^1}(S^2) = \C[x,y]/xy$, and the equivariant parameter $t$ maps to $x+y$. Setting $t=0$ we recover the ordinary cohomology $H^\bullet(S^2)$, which is the endomorphism algebra of $S$ in $\cW(X)$. 

 Of course this observation is not limited to this example, it applies in general to $S^1$-invariant exact compact Lagrangians.
\end{examplex}
 
It seems reasonable to expect that $\cW_{S^1}\!(X)_\lambda$ is generated by invariant Lagrangians, provided that we include the non-compact ones. Moreover, in some examples $\cW(X)$ itself can be generated by invariant Lagrangians; this happens for instance in Example \ref{eg.firsteg} since the cotangent fibre at one of the two fixed points is invariant and generates $\cW(T^*S^2)$. In this situation it may be possible to recover the non-equivariant category $\cW(X)$ from a single spectral component $\cW_{S^1}\!(X)_\lambda$ of the equivariant category.
\pgap

Now suppose we have the set-up of Conjecture \ref{conj.conics}, with (for simplicity) $r=1$. So we have a rank one algebraic torus fibration $\pi: X \to Y$ degenerating over a divisor $D\subset Y$. The conjecture is that:
$$ \cW_{S^1}\!(X)_{-1} \cong \cW(Y\setminus D)$$
Since the category on the left is linear over $\C[t]$, the category on the right should be too. 

There is an obvious guess for what this extra structure on $\cW(Y\setminus D)$ is. Indeed, there is a class $\tau \in SH^2(Y \setminus D)$ corresponding to a simple Reeb orbit going around the divisor $D$ once.  It is sometimes called the Borman-Sheridan class \cite{Tonkonog2}. With this choice of $\tau$, the category $\cW(Y\setminus D)$ becomes linear over $\mathbb{C}[t]$.\footnote{With the same technical caveat about centrality that we raised in Remark \ref{rem.notfibres}.}

\begin{keyconj}\label{conj.tistau} In the situation of Conjecture \ref{conj.conics}, with $r=1$, the action of $t$ on $\cW_{S^1}\!(X)_{-1}$ coincides with the action of $\tau$ on $\cW(Y \setminus D)$. 
\end{keyconj}

If correct this implies that the fibre $\cW(Y\setminus D)|_{t=0}$ agrees with the subcategory of $\cW(X)$ generated by objects $L$ such that $s_L = -\id_L$. If this is the whole of $\cW(X)$ - which looks roughly equivalent to asking that $\cW(X)$ is generated by invariant Lagrangians - then it follows that:
\beq{eq.baseatt=0}\cW(Y\setminus D)|_{t=0}\;\cong\; \cW(X)\eeq
We expect this to hold whenever the wrapped category of $Y$ is zero, which will be true if $Y$ is a subcritical Weinstein manifold such as $\mathbb{C}^n$.  In particular it should work for Example \ref{eg.firsteg}. In the next section we'll discuss this in more detail, including the mirror picture.

When it holds, the equivalence \eqref{eq.baseatt=0} opens a route to proving some new instances of homological mirror symmetry, by bootstrapping up from a theorem about $Y$ to a theorem about $X$.

\begin{rem}
In our discussion after Conjecture \ref{conj.conics} we noted that the deformation
$$ \cW_{S^1}\!(X\setminus \pi^{-1} D)_{-1} \;\leadsto \; \cW_{S^1}\!(X)_{-1} $$
must be trivial. But now we can refine this: the $A_\infty$-structure is not deforming, but the $\C[t]$ structure is.  On the quotient we are starting with the trivial $\C[t]$ structure on $\cW(Y \setminus D)$ and deforming it to the non-trivial one where $t$ acts as $\tau$. The class $\tau$ itself relates to the formal deformation of $\cW(Y \setminus D)$ to the relative category $\cW(Y,D)$, which is the subject of our next section.
\end{rem}

\begin{rem}\label{rem.mirrort}
A more precise version of ``Conjecture'' \ref{conj.mirror} would also include $\C[t]$ structures. If the mirror to $S^1\acts X$ is $\sigma: \mX \to \C^*$ then each category $D^b(\mZ_\lambda)$ is automatically linear over $\C[t]$, because the presentation of $\mZ_\lambda$ as a hypersurface provides a choice of $t\in HH^2(\mZ_\lambda)$. This should agree with the $\C[t]$ structure on $\cW_{S^1}\!(X)_\lambda$. 

If the mirror is an LG model then we again have an automatic choice of 
$$t\in HH^2\big(\MF(\mZ_\lambda, \mW|_{\mZ_\lambda})\big)$$
 but now $t$ might be recording that the superpotential $\mW|_{\mZ_\lambda}$ is varying with $\lambda$, rather than the hypersurface itself. Or they might both be deforming.  
\end{rem}

\subsection{Relative Fukaya categories}\label{sec.relative}

Suppose again we have the setup of Conjecture \ref{conj.conics} with $r=1$.  On the base space $Y$ we can consider the \emph{relative} wrapped Fukaya category:
$$\cW(Y, D)$$
This category is, by construction, linear over a power series ring $\C[[h]]$ where $\deg h = 0$. The relative Fukaya category was first considered by Seidel in \cite{SeidelICM}. The objects are exact Lagrangians in $Y \setminus D$, which may be non-compact but must stay away from $D$, \emph{i.e.}~they should not asymptotically approach to $D$. 
\vspace{3pt}

Setting $h=0$ in the relative category recovers the full subcategory of $\cW(Y\setminus D)$ spanned by these same Lagrangians. If instead we invert $h$, we get the ordinary wrapped Fukaya category $\cW(Y)$, with $h$ essentially becomes the Novikov parameter (\emph{c.f.}~\cite[Lemma 1.2]{LP}). So the relative category captures the `total space' of the formal deformation with central fibre $\cW(Y\setminus D)$ and generic fibres $\cW(Y)$. 
\pgap

Now consider the space $X$. The $S^1$ action makes $\cW(X)$ linear over the ring $\C[s^{\pm 1}]$. At a generic value $\lambda$ of $s$ one expects that $\cW_{S^1}\!(X)_\lambda$ is equivalent to $\cW(Y)$, possibly with a bulk deformation. At $\lambda=-1$, our Conjecture \ref{conj.conics} is that $\cW_{S^1}\!(X)_{-1} \cong \cW(Y\setminus D)$.  
\pgap

Comparing the previous two paragraphs we reach the following conclusion:

\begin{keyconj}\label{conj.relative} 
Suppose we have the setup of Conjecture \ref{conj.conics} with $r=1$. Then the relative wrapped Fukaya category $\cW(Y,D)$ is equivalent to the completion of $\cW(X)$ at $s=-1$. 
\end{keyconj}

If we express this in terms of mirror symmetry as in Section \ref{sec.MSintro} then we are claiming the following. Suppose $X$ is mirror to $(\mX, \mW)$ and the $S^1$ action is mirror to a function $\sigma:\mX \to \C^*$. Then $\cW(Y, D)$ will be equivalent to the category of matrix factorizations on the formal scheme obtained by completing $\mX$ along the divisor $\mZ_{-1}= \sigma^{-1}(-1)$.

\begin{rem} Conjecture \ref{conj.relative} suffers from the same technical defect (1) discussed in Remark \ref{rem.notfibres}: we probably need some $E_2$ structure to be able to complete $\cW(X)$ over $\C[s,s^{-1}]$.  But the more fundamental problem (2) does not apply here, because (unlike when we take fibres) completing the category $\MF(\mX, \mW)$  should give the same result as completing the space $\mX$.
\end{rem}


\begin{rem} The relative category $\cW(Y, \cD)$ is linear over $\C[[h]]$, and the completion of $\cW(X)$ is also linear over a power series ring. We're implicitly claiming that these structures agree, and for this to be compatible with Conjecture \ref{conj.tistau} we need that $h=s+1$ to first-order. But, in general, there might be higher-order corrections. We expect this to be controlled by the quantum Kirwan map (\emph{c.f.}~\cite{WoodwardXu}). 
\end{rem}

\begin{rem}  \label{rem.multiplesingularfibres}One could formulate Conjecture \ref{conj.relative} more generally in the setup of Conjecture \ref{conj.singularfibres}. If there are multiple singular fibres then there will be a relative Fukaya category for each singular value, and the claim is that each is equivalent to the completion of $\cW(X)$ at some corresponding value of $s$. 
\end{rem}


Perhaps the most interesting situation for this conjecture is the case when $\cW(Y)\cong 0$, as we discussed in the previous section. In this case we expect from \eqref{eq.FukayaTeleman} that $\lambda=-1$ is the only non-zero spectral component of $\cW_{S^1}\!(X)$. 

 On the mirror side this is the statement that 
$$\MF\!\big(\mZ_{\lambda}, \mW|_{\mZ_{\lambda}}\big)\cong 0 \quad\quad \mbox{for all } \lambda\neq -1$$
which says that the restriction of $\mW$ to $\mZ_\lambda$ has no critical points if $\lambda \neq -1$. It follows that $\Crit(\mW)$ is contained (set-theoretically) in $\mZ_{-1}$, and hence completing $\mX$ at $\mZ_{-1}$ has no effect on the category of matrix factorizations.\footnote{If the restriction of $\mW$ to some divisor $\mZ$ has no critical points, then it follows that $\mW$ has no critical points in a formal neighbourhood of $\mZ$. But the converse is not true. This elementary fact will be significant in the next section.} Equivalently, the subcategory $\MF_{\mZ_{-1}}(\mX, \mW)$ of matrix factorizations supported on $\mZ_{-1}$ is in fact the whole of $\MF(\mX, \mW)$. 

Returning to the A-side, we conclude that if $\cW(Y)\cong 0$ then completing $\cW(X)$ at $s=-1$ has no effect. Equivalently, $\cW(X)$ is generated by objects $L$ such that $s_L= -\id_L$. So \eqref{eq.baseatt=0} holds, and Conjecture \ref{conj.relative} becomes the claim that $\cW(Y,D)$ and $\cW(X)$ are equivalent.
\pgap 


To understand this claim in more detail we draw the following diagram:
\begin{equation}\label{adjunctiondiagram}
\begin{tikzcd}
\MF(\mX, \mW) \ar[shift left=0.5ex]{d}{i^*} \ar[dashrightarrow, dash]{r}{\simeq}  &  \cW(X) \ar[r, "\simeq", "F"']
 \ar[dashrightarrow, shift left=0.5ex]{d}& \cW(Y, D) \ar[shift left=0.5ex]{d}{h=0}   \\
\MF(\mZ_{-1}, \mW|_{\mZ_{-1}}) \ar[shift left=0.5ex]{u}{i_*}  \ar[dashrightarrow, dash]{r}{\simeq} & \cW_{S^1}\!(X)_{-1} \ar[shift left=0.5ex]{u}{t=0}& \cW(Y\setminus D) \ar[l, "G", "\simeq"'] \ar[dashrightarrow, shift left=0.5ex]{u} 
\end{tikzcd}
\end{equation}

Here the categories along the top row are all supposed to be equivalent, as are the three categories along the bottom. The vertical arrows are an adjunction. 

On the B-side the adjunction is obvious, it's given by pushing forward or pulling back along the inclusion $i: \mZ_{-1}\into \mX$. On the symplectic manifold $X$ only the upward functor is obvious - it's the functor from Section \ref{sec.recovering} which takes an invariant Lagrangian in $\cW_{S^1}(X)_{-1}$ and regards it as an object of $\cW(X)$ instead. The adjoint to this functor is less clear. 

Curiously, on the base $Y$ it is the downward functor which is obvious - just set the deformation parameter $h$ to zero - and the upward functor which is a little more surprising.

Our conjectured equivalence $G: \cW(Y\setminus D) \isoto \cW_{S^1}\!(X)_{-1}$ should roughly be given by the Lagrangian correspondence $\Gamma$ associated to $\mu^{-1}(0)$, see \eqref{momentLag}. So, at least approximately, it sends a Lagrangian $\Lambda$ to the invariant Lagrangian $\pi^{-1}(\Lambda)\cap \mu^{-1}(0)$.

On the top line, the functor $F:  \cW(X) \isoto \cW(Y, D)$ needs to be adjoint to $G$, so it's given by the same correspondence $\Gamma$ but viewed as a functor in the opposite direction. To apply this functor to a Lagrangian $L\subset X$ we must perturb until $L$ intersects $\mu^{-1}(0)$ transversely, then we project $L\cap \mu^{-1}(0)$ down to $Y$. 

\begin{rem} 
If we drop the condition that $s=-1$ is the only non-zero component then we can still draw this diagram, the only change is that $F$ is no longer an equivalence. It is mirror to the restriction from $\mX$ to the formal completion of $\mX$ along $\mZ_{-1}$. 
\end{rem}

\begin{egu} Revisiting Example \ref{eg.firsteg}, Conjecture \ref{conj.relative} claims that we have an equivalence:
\beq{eq.relativeeq} \cW(T^*S^2) \cong \cW(\C, \{\pm 1\})\eeq
In this example we can compute the effects of all functors in \eqref{adjunctiondiagram} quite explicitly.

Take a straight arc in $\C$ from either puncture to infinity, \emph{i.e.} either $\Lambda_2$ or $\Lambda_3$ in Figure \ref{fig.arcs}. Under $G$ they lift to the cotangent fibres at the two fixed points. These objects are not isomorphic in $\cW_{S^1}\!(X)_{-1}$, but they become isomorphic when we pass to $\cW(X)$ since all cotangent fibres are isomorphic there. If we want to now apply the functor $F$ we must take a cotangent fibre at a generic point so that it intersects $\mu^{-1}(0)$ transversely, and the result is the arc $\Lambda_4$. So the `surprising' functor $\cW(Y\setminus D)\to \cW(Y, D)$ appearing in the third column of \eqref{adjunctiondiagram} sends both $\Lambda_2$ and $\Lambda_3$ to $\Lambda_4$. 

As a test of our claimed equivalence \eqref{eq.relativeeq} we can compute  $\Hom_{\mathcal{W}(Y,D)}(\Lambda_4,\Lambda_4)$ directly using the wrapped Floer complex. The result is 
\[ \frac{\mathbb{C}\langle u,v\rangle[h]}{(u^2,v^2)}, \ \ \  du = h\cdot 1, dv = h \cdot 1, \ \ |u|=|v|=-1 \]
and it is easy to see that this is quasi-isomorphic to $\mathbb{C}[x]$ with $|x|=-1$, which is indeed the endomorphism algebra of a cotangent fiber in $\mathcal{W}(T^*S^2)$.

If we apply $G$ to the arc $\Lambda$ from Example \ref{eg.firsteg} then we get $S^2$, as an object of $\cW_{S^1}\!(X)_{-1}$. But applying $F$ to $S^2 \in \cW(X)$ does not produce $\Lambda$;  no arc ending at $1$ or $-1$ can be in the image of $F$.  In fact $F(S^2)$ is a figure-eight encircling both punctures.

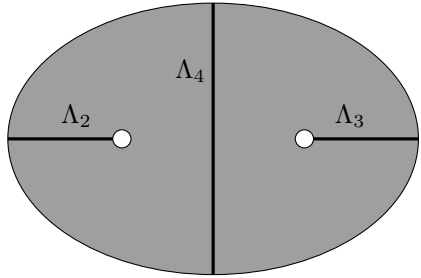
\begin{figure}[h!]
\label{fig.arcs}
  \begin{tikzpicture}[scale=0.6]

\draw[fill=gray!75, x radius=4.5, y radius=3] (4.5, -3) ellipse;

\draw[ fill=white,radius=.2] (2.5,-3) circle ;

\draw[ fill=white,radius=.2] (6.5,-3) circle ;

\draw[very thick] (0,-3)-- (2.3,-3);

\draw[very thick] (6.7,-3) -- (9,-3);

\draw[very thick] (4.5,-6) -- (4.5,0);

\node at (1.5,-2.5) {$\Lambda_2$};

\node at (7.5,-2.5) {$\Lambda_3$};

\node at (4,-1.5) {$\Lambda_4$};
\end{tikzpicture}
\caption{Lagrangians in $\mathbb{C} \setminus \{\pm 1\}$.}
\end{figure}

\end{egu}

\begin{rem} The `surprising' functor $\cW(Y\setminus D) \to \cW(Y, D)$ formally is given by setting the Borman-Sheridan class $\tau$ to zero. We can get a more explicit understanding by passing through the mirror. 

 On the B-side, the composition $i_*i_*$ acts on objects as
$$i^*i_*\!:\; \cE \;\mapsto\; \big[ \cE \stackrel{t_\mathcal{E}}{\To} \cE[2]  \big]$$
where $t_\mathcal{E}$ is the $HH^2$ class from Remark \ref{rem.mirrort} evaluated on the object $\cE$.  On the A-side, this corresponds to taking a Lagrangian $L \subset Y\setminus D$ and forming the twisted complex
$$ \big[ L \stackrel{\tau_L}{\To} L[2]  \big] $$
where $\tau_L$ is the image of $\tau$ under the closed-open map $\mathcal{CO}^0 : SH^*(\mathcal{W}(Y \setminus D)) \to HF(L,L)$. But we claim that this twisted complex can be represented by an object $\tilde{L}$ which does not asymptote to $D$, so it can be viewed as object of the relative category $\cW(Y, D)$.  On the level of objects we believe this is our functor. It's easy to see that this consistent with our discussion of the objects $\Lambda, \Lambda_2$ and $\Lambda_3$ in the example above. 

Note that if $L$ does not have an end converging to $D$ then we should have $\tau_L=0$, and so the result of applying the functor should be just $L \oplus L[1]$. On the mirror this is what happens when $\cE=i^*\cF$ for some $\cF\in \MF(\mX, \mW)$. 
\end{rem}

\begin{rem} \label{rem.relativeeq} We could also consider the relative category:
 $$\cW(X, \pi^{-1} D)$$
 This will be linear over $\C[[h]]$, and it should also have a $\C[s, s^{-1}]$ structure coming from the $S^1$ action. We believe that these two structures are essentially the same, and this encapsulates many of our claims. See Example \ref{eg.simplest2} for an example (on the mirror side) of what we mean by this.
\end{rem}

\subsection{Topological $S^1$ actions on categories and the KRS 2-category}\label{sec.topologicalactions}

The structures that we have encountered can be understood more clearly using the framework of `topological group actions' on categories sketched by Teleman in \cite{Teleman}. In this section we'll give a brief explanation of these ideas, restricting only to the case of topological $S^1$ actions. The generalization to actions of higher rank tori is easy; the non-abelian case is considerably deeper but we won't need it. 

To first approximation a topological $S^1$ action on a category $\cC$ is a choice of natural automorphism $s: \id_\cC \to \id_\cC$ of the identity functor, or equivalently a $\C[s, s^{-1}]$ structure on $\cC$. We saw in Section \ref{sec.MSintro} that the Seidel element provides this structure on $\cC=\cW(X)$ when $X$ is a Hamiltonian $S^1$-manifold. On the mirror side, a non-vanishing function on $\mX$ provides this structure on $\cC=D^b(\mX)$. 

 A heuristic general explanation for this definition goes as follows. Imagine that $\cC$ is a topological category (or $\infty$-category) with a space of objects $\Ob_{\cC}$. Then the space of automorphisms of any object $x\in \Ob_{\cC}$ is, by definition, the based loop space $\Omega_x\Ob_{\cC}$. A `topological $S^1$ action' means an action $S^1 \acts \Ob_{\cC}$, and then for any object $x$ we have an automorphism $s_x:x\isoto x$ given by the orbit of $x$ under the $S^1$ action.

If we have a topological $S^1$ action on $\cC$ we can then ask about the equivariant category. A natural definition of an `equivariant structure' on an object $x$ is a choice of contracting homotopy $h$ of the loop $s_x$ (this is analogous to the definition of an equivariant sheaf). Such an $h$ induces - up to homotopy - an $S^1$ action on $\Omega_x\Ob_{\cC}$, and then the automorphism space of $(x,h)$ in the equivariant category will be something like the homotopy quotient $\Omega_x\Ob_{\cC}\times_{S^1} ES^1$. Hence, if we linearize by replacing spaces with their cohomology, we expect the equivariant category to be linear over $H^\bullet(BS^1)=\C[t]$. 
 
More algebraically, this discussion suggests that the equivariant category should be the fibre of $\cC$ at the point $s=1$. And indeed, the fibre category is always linear over $\C[t]$.\footnote{This is for the same reason that the derived category of a hypersurface is linear over $\C[t]$, see Remark \ref{rem.mirrort}.} More generally we could take the fibre at some value $s=\lambda$ and interpret this as a `twisted equivariant category'. This looks like it could be the correct framework for the equivariant Fukaya categories $\cW_{S^1}\!(X)_\lambda$ and the mirror symmetry picture from Section \ref{sec.MSintro}.

Unfortunately (as we warned in Remark \ref{rem.notfibres}) this story is too simplistic.
 
\begin{eg}\cite[Example 4.5.]{Teleman}
Let $S^1$ act on $\C^n$ with the standard diagonal action, so the non-empty symplectic quotient is $\P^{n-1}$. Then $\Fuk(\P^{n-1})$ is a spectral component of $\cW(\C^n)$. But $\cW(\C^n)\cong 0$ and $\Fuk(\P^{n-1})\ncong 0$, so the latter cannot possibly be a fibre of the former.
\end{eg}

That example cannot be $\Z$-graded, but the next one can.

\begin{eg}\label{eg.simplest1}
The simplest example for Conjecture \ref{conj.conics} is $X=\C^2$ mapping to $Y=\C$ via the map $\pi=xy$, as we considered in Example \ref{eg.simplest0}. Our conjecture is that $\cW_{S^1}\!(X)_{-1} \cong \cW(\C^*)$. Since this is non-zero, it cannot be a fibre of $\cW(X)$. 
\end{eg}

On the mirror side these two examples are quite similar: the mirror superpotential $\mW$ has no critical points, but the restriction of $\mW$ to a hypersurface $\mZ_\lambda$ does have critical points.\footnote{In the second example, $\mW|_{\mZ_{-1}}$ is identically zero, see Example \ref{eg.simplest2}.} So $\MF(\mW)\cong 0$ but $\MF(\mW|_{\mZ_\lambda})\ncong 0$. This is obviously a general phenomenon and illustrates sharply why $\MF(\mW|_{\mZ_\lambda})$ cannot be a fibre of $\MF(\mW)$. 

\newcommand\KRS{K\!RS}

Teleman's solution to this puzzle uses the \textit{Kapustin-Rozansky-Saulina 2-category}. This is a conjectural construction from 3D topological field theory (Rozansky-Witten theory) that associates a 2-category to a holomorphic symplectic manifold $M$.

 An object of this 2-category $\KRS(M)$ is a choice of holomorphic Lagrangian $L\subset M$, together with a category $\cC_L$ which is linear over $L$, \emph{i.e.}~which carries an action of the tensor category $Coh(L)$. In particular each $L$ defines an object by setting $\cC_L = Coh(L)$. These special objects should generate $\KRS(M)$, so using the Yoneda embedding we can view a general object $\mathfrak{X}$ as an operation 
$$L \mapsto \Hom_{\KRS(M)}(L, \mathfrak{X})$$
which associates a (1-)category to each Lagrangian. We need this perspective for our version of equivariant mirror symmetry (see below). It also has the advantage that we can immediately handle objects $\mathfrak{X}$ supported on \emph{singular} Lagrangians, by viewing them as on operation on smooth Lagrangians. 

For topological $S^1$ actions the relevant holomorphic symplectic manifold is the shifted cotangent bundle: 
$$M = T^\vee[-2]\mathbb{G}_m = \Spec \C[s^{\pm 1}, t] $$
In the $\Z_2$-graded world this is symplectic, in the $\Z$-graded world it is shifted-symplectic. The latter is much easier because then there are very few (smooth) Lagrangians. Indeed the only possible $L$'s are the zero section and the cotangent fibres, and these intersect in a very simple way. So in the $\Z$-graded world, an object of $\KRS(T^\vee[-2]\mathbb{G}_m)$ should be the following data:
\begin{enumerate}\setlength{\itemsep}{2pt}

\item[(i)] A $\C[s^{\pm1}]$-linear category $\cC$.
\item[(ii)] A $\C[t]$-linear category $\cD_\lambda$ for each $\lambda\in \C^*$. 
\item[(iii)] For each $\lambda$ an embedding
$$\cC|_\lambda  \;\into\; \cD_\lambda$$
where the source is the fibre of $\cC$ at $s=\lambda$, and the image is the subcategory of $\cD_\lambda$ consisting of objects supported at $t=0$. 
\item[(iv)] For each $\lambda$ an embedding
$$(\cD_\lambda)|_0 \;\into\; \cC$$
where the source is the fibre of $\cD_\lambda$ at $t=0$ and the image is the subcategory of $\cC$ consisting of objects supported at $s=\lambda$. 
\end{enumerate}
\begin{rem}(iii) and (iv) might be equivalent. They arise from considering the morphisms between the zero section and the cotangent fibres in $\KRS(T^\vee[-2]\mathbb{G}_m)$, which we have not discussed. \end{rem}

One of Teleman's insights is that a ``category with a topological $S^1$-action'' is really an object of $\KRS(T^\vee[-2]\mathbb{G}_m)$. In the $\Z$-graded world we can view it as a system of categories as described above.

 On the B-side we get an example by taking a smooth variety $X$ with a function $\sigma: X \to \C^*$ and setting 
$$\cC=D^b(X) \aand \cD_\lambda = D^b(Z_\lambda)$$
where $Z_\lambda = \sigma^{-1}(\lambda)$. Here (iii) becomes the fact that the fibre of $D^b(X)$ at $s=\lambda$ is $\Perf(Z_\lambda)$, which is also the subcategory of $D^b(Z_{\lambda})$ on which $t$ acts nilpotently. And (iv) is the fact that setting $t=0$ in $D^b(Z_\lambda)$ recovers the subcategory $D^b_{Z_\lambda}\!(X)\subset D^b(X)$ of objects supported on $Z_\lambda$.  

More generally, we can take an LG-model $(X, W)$, with R-charge, equipped with a function $\sigma: X \to \C^*$. Then we set $\cC=\MF(X, W)$ and $\cD_\lambda = \MF(Z_\lambda, W|_{Z_\lambda})$. 

On the A-side, we should get an example by taking a graded Hamiltonian $S^1$-manifold $X$ and setting:
$$\cC = \cW(X)\aand \cD_\lambda=\cW_{S^1}\!(X)_\lambda$$
 If all fibres of the moment map are regular we should be able to replace each $\cW_{S^1}\!(X)_\lambda$ by $$\cW(X\sslash_{\log |\lambda|} S^1, \,\arg (\pm\lambda))$$ as discussed in Section \ref{sec.Hamiltonianreduction}.  In the situation of Conjecture \ref{conj.conics} we should get an example if we set $\cC= \cW(X)$, set $\cD_\lambda = \cW(Y)$ for $\lambda\neq -1$, and set $\cD_{-1}= \cW(Y\setminus D)$. 

In fact on the A-side we can be more direct. The full equivariant category $\cW_{S^1}(X)$ should in a natural way be an object of $\KRS(T^\vee[-2]\mathbb{G}_m)$, supported on some Lagrangian $L$, and it should induce the systems of categories above. But we don't know the analogue of this on the B-side, \emph{c.f.} Remark \ref{rem.equivariantmirror}.

\begin{remu} If we are $\Z$-graded, and each category $\cD_\lambda$ is smooth and proper, then $t$ must be torsion (since $\deg(t)=2$). It follows that (iii) must be an equivalence $\cC|_\lambda \isoto \cD_\lambda$. In this situation the category $\cC$ determines everything and the simpler story described earlier in this section is enough, \textit{c.f.}~\cite[Remark 3.2]{Teleman}. In the KRS 2-category we just have the category $\cC$, supported on the Lagrangian $\{t=0\}$. 

 But Example \ref{eg.simplest1} demonstrates that the simpler story fails under just the assumption of $\Z$-grading alone. In that example we have $\cD_{-1} = \cW(\C^*)$, and we claim that the $\C[t]$ structure is just the co-ordinate function, \emph{i.e.}
$$\cD_{-1} \cong \Spec \C[t^{\pm 1}]$$
(see Example \ref{eg.simplest2}). So $\cD_{-1}\ncong 0$, but it contains no objects supported at $t=0$. This object of $\KRS(T^\vee[-2]\mathbb{G}_m)$ is supported on the incomplete Lagrangian $\{s=-1, t\neq 0\}$. 
  
\end{remu}

\section{Toric mirror symmetry}\label{sec.toricms}

In this section we first review some heuristics in mirror symmetry, which are well-known, and possibly all established.\footnote{But chasing the references throughout the vast literature on this subject is beyond the capabilities of the authors.} We then discuss a series of examples giving evidence for the conjectures of Section 1.  
\pgap

\subsection{General heuristics}\label{sec.generalheuristics}

The starting point is a torus $T\cong (\C^*)^n$. The mirror to $T$ is just the dual torus $T^\vee$. From here we get to more complex examples by applying the following two constructions:

\begin{enumerate}[leftmargin=0.6cm, itemsep=10pt]
\item \label{h1}
We can partially compactify $T$ to a toric variety $X$. On the mirror side this corresponds to adding a superpotential $\mW \in \Gamma(\cO_{T^\vee})$, so the mirror to $X$ is a Landau-Ginzburg model $(T^\vee, \mW)$.  More precisely, the rays of the toric fan for $X$ give a finite set of points in the lattice
$$ L = \Hom(\C^*, T) = \Hom(T^\vee, \C^*) $$
and $\mW$ is the corresponding Laurent polynomial.\footnote{A more sophisticated version includes toric stacks.}

If we start on the A-side with the wrapped Fukaya category $\cW(T)$ then adding in the toric boundary $\partial X$ deforms the  category, and this corresponds on the mirror to deforming $D^b(T^\vee)$ to the category of matrix factorizations $\MF(T^\vee, \mW)$. 

In the other direction, adding $\mW$ to $T^\vee$ changes $\cW(T^\vee)$ by introducing stops at infinity. This means that certain non-compact Lagrangians become non-zero objects, and these correspond to complexes in $D^b(X)$ supported on $\partial X$. 

One can of course swap the roles of $T$ and $T^\vee$ here, and we can also do both operations to both sides, so the mirror to a  toric LG model $(X, W)$ is another toric LG model $(\mX, \mW)$.

\item \label{h2} Suppose $X$ is the total space of a line bundle $L$ over a base $B$. Given a section $f\in \Gamma(B, L^\vee)$ of the dual line bundle we have an induced superpotential $W = pf$ on $X$, where $p$ is the fibre co-ordinate. Then the LG model $(X,W)$ is equivalent to the variety:
$$ H = \{f=0\}\subset B $$
More generally one can introduce a function $g\in \Gamma (\cO_B)$ and set $W'=pf + g$, then the LG models $(X, W')$ and $(H, g|_H)$ are equivalent.

On the B-side this phenomenon is called Kn\"orrer periodicity, and it's been proven by many people (\cite{OrlovKP},...) that $\MF(X, W')\cong \MF(H, g|_H)$. On the A-side the analogous result has also been established recently \cite[Thm.~1]{Jeffs}.

\end{enumerate}

To these two well-established heuristics we add a third, following Section \ref{sec.MSintro}.

\begin{enumerate}[leftmargin=0.6cm, itemsep=10pt]\addtocounter{enumi}{2}

\item \label{h3}Suppose we have a mirror pair of toric LG models $(X, W), (\mX, \mW)$ as in (1). A choice of 1-parameter subgroup $l\in L$ determines an $S^1$ action on $X$, and also a monomial function $l\in \Gamma(\cO_{\mX})$.  The mirror to the $S^1$ action is the function:
$$\sigma = l $$
We require that $W$ is invariant under the $S^1$ action, meaning that every monomial in $W$ lies in the sublattice $\langle l\rangle^\perp \subset L^\vee$. It follows that $\sigma$ is a non-vanishing function on $\mX$. 
\end{enumerate}

\begin{rem} These heuristics only work as stated under the assumption that $X$ is exact or monotone, or more generally satisfies Assumption 1.4 of \cite{AAK}, which guarantees that $\mW$ and $\sigma$ don't receive any higher-order corrections involving Novikov parameters.
\end{rem}

\begin{rem} Arguably our presentation here is backwards. Since a toric variety is a symplectic reduction of a vector space, you can derive heuristic \eqref{h1} from \eqref{eq.FukayaTeleman}, ``Conjecture'' \ref{conj.mirror}, heuristic \eqref{h3}, and the knowledge that the mirror of $\C^N$ is $((\C^*)^N\!,\, \mx_1+...+\mx_N)$. And this seems closer to the original physical derivation of toric mirror symmetry. 
\end{rem}

\begin{eg} \label{eg.simplest2} Recall the very simple Example \ref{eg.simplest0} where $X=\C^2$. We view it as fibred over $Y=\C$ via the map $\pi =xy$ which has a single singular fibre over $D=\{0\}$. The fibration is preserved by the $S^1$ action on $X$ with weights $(1,-1)$. 

Following heuristic (\ref{h1}) above, the mirror to $X$ is the LG model $\mX=(\C^*)^2$ with $\mW = \mx + \my$. And heuristic (\ref{h3}) says that the mirror to the $S^1$ action is the function $\sigma = \mx/\my$. If we want to have a $\Z$-grading the only option is that $\deg(\mx)=\deg(\my)=2$, and then $\deg(\mW)=2$ but $\deg(\sigma)=0$, as required.

 If we want to make the $\C[s^{\pm 1}]$ structure completely explicit we can change co-ordinates and declare that the mirror to $X$ is 
\beq{eq.mirror1} (\C^*_{\my} \times \C^*_s, \, (1+s)\my) \eeq
 If we slice at $s=-1$ we just get $\C^*_{\my}$ with no superpotential. So our ``Conjecture'' \ref{conj.mirror} predicts:
$$\cW_{\cS^1}\!(X)_{-1} \cong D^b(\C^*)$$
 This fits Conjecture \ref{conj.conics} since $Y\setminus D = \C^*$ is self-mirror. If we slice at value $s\neq -1$ we get the mirror to $Y=\C$, and indeed the category is zero on either side of the mirror. This supports the claims we made in Example \ref{eg.simplest0}.

The $\C[t]$ structure on $D^b(\C^*_{\my})$ is evidently just the co-ordinate function, $t=\my$ (see Remark \ref{rem.mirrort}). Note that this is consistent with the grading since $\deg(\my)=2$. It also matches the Borman-Sheridan class on $\cW(\C^*)$, which counts the additional holomorphic disc that appears once we fill in the puncture. From the same picture we can see that the relative Fukaya category on the base must be
$$\cW(Y, D) \cong \MF\!\big(\Spec \C[\my^{\pm 1}][[h]], \; h\my\big) $$
which is indeed the completion of \eqref{eq.mirror1} at $s=-1$, matching Conjecture \ref{conj.relative}. But this last observation is vacuous, since both categories are zero.

We can look more carefully at our proposed equivalence $\cW(\C^*) \isoto \cW_{S^1}\!(X)_{-1}$, \emph{i.e.} the functor $G$ from \eqref{adjunctiondiagram}. If $\Lambda\subset \C^*$ is an arc from 0 to $\infty$ then it's clear that $G(\Lambda)$ will be an $S^1$-invariant Lagrangian thimble $L$. The claim is that such an $L$ is a generator of $\cW_{S^1}\!(X)_{-1}$ and has endomorphism algebra $\C[t^{\pm 1}]$. Presumably one could check this claim directly. 

Finally, let us fulfil our promise in Remark \ref{rem.relativeeq} and discuss the relative category $\cW(X, \pi^{-1}D)$. 
Since $\pi^{-1}(D)$ here is just the toric boundary $\{\mx\my=0\}$, it's clear that the relative category is mirror to:
$$ \MF\!\Big(\C[\mx^{\pm 1}, \my^{\pm 1}][[h]], \,h(\mx + \my)\Big) \; \; \cong\; \;  \MF\!\Big(\C[\my^{\pm 1}, s^{\pm 1}][[h]], \,h(1+s)\my\Big)$$
Here we have applied the same co-ordinate change as above.  We see that the $\C[s^{\pm 1}]$ structure  and the $\C[[h]]$ structure are tied together in that they act through the product $h(s+1)$. It follows that once we slice at $s=-1$ the deformation in the $h$ direction is trivial; this is our claim that adding in $\pi^{-1}(D)$ does not deform the equivariant Fukaya category of $X\setminus D$. But the $\C[t]$ structure on the slice does deform, because $t=h\my$. This should be a general phenomenon.

\end{eg}

\subsection{Torus fibrations over a torus}

\begin{egu} \label{eg.nodedeforming2}For a variation on Example \ref{eg.simplest2}, set $X$ to be the open variety $X= \C^2\setminus \{ xy=1\}$ and consider the projection:
$$\pi=xy-1\! :\; X \to Y=\C^*$$
We considered exactly this $X$ in Example \ref{eg.nodedeforming1}. It's an algebraic torus fibration with a single singular fibre over $D=\{1\}$ and it's again preserved by the $S^1$ action with weight $(1,-1)$.  The difference with Example \ref{eg.simplest2} is that now $Y=\C^*$ is log-Calabi-Yau and $\cW(Y)\neq 0$. So we expect that $\cW(X)$ will be supported over the whole of $\C^*_s$ and not just at $s=-1$. 

To apply our mirror symmetry heuristics we view $X$ as the hypersurface:
$$X = (z-xy+1)\;\subset \C^2_{x,y}\times \C^*_z$$
 Then by heuristic (\ref{h2}) it is equivalent to the LG model:
$$ \C^3_{x,y,p} \times \C^*_z, \quad W = p(z-xy + 1) $$
Without $W$ the mirror would be $(\C^*)^4$ with superpotential $\mW = \mx + \my + \mp$, and $\sigma = \mx/\my$. The presence of $W$ induces a partial compactification of the mirror. To see what this is we perform the co-ordinate change
$$ z \mapsto pz, \quad y \mapsto pxy $$
so now $W = z - y + p$. On the mirror we have to perform the dual co-ordinate change:
$$ \mp \mapsto \mp\mz\my,  \quad \mx \mapsto \mx\my$$   
So the mirror is 
\beq{nodesmoothing}\C^3_{\my, \mz, \mp}\times \C^*_{\mx}, \quad \mW = (1+\mx)\my + \my\mz\mp\eeq
with $\sigma = \mx$. 



Using heuristic \eqref{h3} we can remove the $\my$ direction and get that the mirror is the hypersurface
$$ \mH =  (1+\mx+\mz\mp)\; \subset \C^2_{\mz,\mp}\times \C^*_{\mx} $$
which is in fact isomorphic to $X$. If we view $\mH$ as an open subset of $\C^2_{\mz, \mp}$ then $\sigma = -1 -\mz\mp$, as claimed in Example \ref{eg.nodedeforming1}. We saw there how slicing at different values of $\sigma$ supports our conjectures.

Completing $\mH$ at $\lambda=-1$ produces the formal smoothing of the node, and this agrees with the relative Fukaya category $\cW(Y, D)$ as in Conjecture \ref{conj.relative}. But it is not the same as $\mH$, reflecting the mirror fact that $\cW(X)$ is not generated by invariant Lagrangians. The subcategory generated by invariant Lagrangians, which is the fibre of $\cW_{S^1}\!(X)_1\cong D^b(\mz\mp=0)$ at $t=0$, should be equivalent to the subcategory 
 $$D^b_{\mz\mp=0}(\mH) \subset D^b(\mH)$$
 of objects supported on the node. 
\end{egu}

The previous example generalizes immediately to rank 1 algebraic torus fibrations over a torus $Y= (\C^*)^n$. Given a divisor $D=(f)\subset Y$ we have an algebraic torus bundle
$$X = (xy - f) \subset \C^2_{x,y}\times Y  $$
with an $S^1$ action having weights $(1,-1)$ on $(x,y)$ and fixing $Y$. 
\pgap

First we construct the mirror to $D$. Since $D$ is a hypersurface in the torus $Y$ we can use heuristic \eqref{h2} to replace it with the LG model $(Y\times \C_p, pf)$. If we consider instead $(Y\times \C^*_p, pf)$ then the mirror would be a toric Calabi-Yau $\mD$ of dimension $n+1$, determined by the Newton polytope of $f$. Adding the divisor $p=0$ induces a superpotential $\mW$ on the mirror, and $\mW$ is the monomial cutting out the toric boundary $\partial \mD$.
\pgap

Next we consider the mirror to $Y\setminus D$. This is the hypersurface $(z-f)\subset \C^*_z\times Y$ so it's equivalent to the LG model $W = p(z-f)$ on $\C_p\times \C^*_z\times Y$. Then the mirror is $\mD \times \A^1_{\mz}$ with the superpotential $\mW\mz$. By heuristic \eqref{h2} we deduce that the mirror to $Y\setminus D$ is the toric boundary of $\mD$. 
\pgap

Finally we can discuss the mirror for $X$. Once again we replace $X$ with an LG model $\C^3_{x,y,p}\times Y$ with $W = (xy-f)p$. The mirror is $\mD \times \C_{\my}\times \C^*_{\mx}$ with superpotential
$$ \mx\my + \my + \my\mW $$
and the mirror to the $S^1$ action is $\sigma = \mx$. Using heuristic \eqref{h2} and writing $s=\mx$ we get that the mirror to $X$ is the hypersurface:
$$\mX =  (\mW + 1 + s ) \; \subset \mD \times \C^*_s $$
This is isomorphic to $\mD\setminus \{\mW=-1\}$. Projecting onto $s$ exhibits it as a family of tori, degenerating at $s=-1$ to the toric boundary $\partial \mD$. This matches our conjectures that a generic fibre of $\Fuk(X)$ should be $\Fuk(Y)$ and the fibre at $s=-1$ should be $\cW(Y\setminus D)$. 
\pgap

If we want the relative Fukaya category $\Fuk(Y,D)$ we identify $Y\setminus D$ with the hypersurface $(z-f)$ as above, and note that $D$ is now the divisor $z=0$. So the mirror to $\Fuk(Y,D)$ is
$$\mD \times \A^1_{\mz} \times \operatorname{Spf} \C[[h]], \quad  \mW\mz + h\mz $$
which is equivalent to the completion of $\mX$ at $s=-1$. 

\subsection{Torus fibrations over affine space}

The construction of the previous section can be adapted to rank 1 algebraic torus fibrations over affine space by partially-compactifying $Y=(\C^*)^n$ to $\overline{Y}=\C^n$. On the mirror side this introduces $n$ extra monomials in the superpotential. 

\begin{eg}
In Example \ref{eg.nodedeforming2}, we can extend $Y$ to $\overline{Y}=\C$ and $X$ to $\overline{X}=\C^2$. Obviously this just a translation of Example \ref{eg.simplest2}, but in our mirror description this is slightly less obvious. Adding in the divisor $\{z=0\}$ adds a monomial to $\mW$, and the mirror becomes:
$$\C^3_{\my, \mz, \mp}\times \C^*_{\mx}, \quad W = \mx\my + \my + \my\mz\mp + \mz  $$
Using heuristic \eqref{h2} this is equivalent to the hypersurface
$$(\my\mp + 1) \subset \C^2_{\my,\mp}\times \C^*_{\mx}$$
equipped with the superpotential $\mx\my + \my$. This is the mirror from Example \ref{eg.simplest2}. 
\end{eg}

\begin{egu}\label{eg.TS2again}
We started this paper with Example \ref{eg.firsteg}, the symplectic manifold $X=T^*S^2$ obtained as the hypersurface $(xy -z^2 + 1) \subset \C^3$. This is an algebraic torus fibration over $Y=\C$ degenerating over the two points $D=(z^2-1)$. Now let's find the toric mirror. 

We start with the mirror to $D$, viewed as a divisor inside $\C^*$. Following the previous section the mirror is the LG model
$$\mD = \operatorname{Tot}\big\{ \cO(-2) \to \P^1_{\um:\mv}\big\},\; \; \mW = \um\mv\mp $$
where $\mp$ is the fibre co-ordinate. Note that $\mW$ has two isolated non-degenerate critical points, which is consistent with it being mirror to $\mD$.  
 
The mirror to $\C^*\setminus D$ is the toric boundary $\partial \mD$. The mirrors to punctured surfaces are well-studied, and $\partial \mD$ is the `balloon-chain' mirror, see for example \cite{LPorder}.

Let $X^o$ be the intersection of $X$ with the torus, so it's an algebraic torus fibration over $\C^*$, and it's $X$ with one smooth fibre cut out. The mirror to $X^o$, still following the previous section, is:
\beq{eq.mirror2}\mX^o =\mD \setminus \{\um\mv\mp = -1\}\eeq
 So $X^o$ is self-mirror. 
 
To get the mirror to $X$ itself we take the LG model mirror to $X^o$ and add one extra term to the superpotential, dual to $\{z=0\}$.  The result is:
$$\mD \times \C_{\my}\times \C^*_{s}, \quad \my(1+s) + \my\um\mv\mp + \mp\um^2$$
It's a little harder to interpret this model geometrically but we can make some observations:

\begin{itemize}\setlength{\itemsep}{5pt}

\item If we fix the value of $s$ to any value except $s=-1$ then the superpotential on that divisor has no critical points, so the category of matrix factorizations is zero. This fits our conjectures since $\Fuk(Y)\cong 0$. 

\item On the divisor $\{s=-1\}$ we can use heuristic \eqref{h2} to remove the $\mp$ direction and we get the hypersurface
$$ \mH = (\my\um\mv + \um^2) \subset \P^1_{\um:\mv}\times \C_{\my} $$
which is an affine variety isomorphic to the node. This fits our conjectures, since $Y\setminus D$ is the pair-of-pants. 

\item On any slice $\{s=\lambda\}$ the critical points of the superpotential of are contained in the affine chart where $\mv\neq 0$. Therefore we can delete the locus $\{\mv=0\}$ and get an equivalent model. After an obvious co-ordinate change the result is: 
\beq{nodesmoothing2}\C^3_{\um, \mp, \my}\times \C^*_s, \quad (\my-\mp)(1+s) + \my\um\mp\eeq

\item It's instructive to compare this with Example \ref{eg.nodedeforming2}, specifically the LG model \eqref{nodesmoothing}. There we saw the derived category of the node appearing at $s=-1$, and the generic fibre was a smoothing of the node. Here we have a more abstract deformation of the node where the generic fibre is zero. 

 On the A-side one can deform the wrapped category of the pair-of-pants either by capping off one puncture or by capping off two punctures simultaneously. Example \ref{eg.nodedeforming2} was the former, this is the latter.


\end{itemize}

\end{egu}

In all such examples we have $Y= \C^n$ and hence $\cW(Y)\cong 0$. As we explained in Section \ref{sec.relative}, it should follow that
$$\Crit(\mW)\subset \{s=-1\}$$
and it easy to check that this is indeed true in the two examples above.

\subsection{Hypertoric quotients}

\begin{egu}\label{eg.hypertoric1} The manifold $T^* S^2$ is not just symplectic but actually hyperk\"ahler, and can be obtained as a hyperk\"ahler quotient of $\C^4$. This gives another perspective on Example \ref{eg.TS2again}.

Take the hypersurface
 $$X= (xy + zw  - 1) \subset \C^4 $$
which symplectically is $T^*S^3$. If we let $U(1)$ act on $\C^4$ with weights $(1,-1,1,-1)$ then $X$ is preserved and the action on it is free. Taking the symplectic quotient at the moment map value $\mma=0$ gives $T^*S^2$. At other values of $\mma$ we get a non-exact symplectic manifold which is diffeomorphic to $T^*S^2$, and it is known that the wrapped category of all these non-exact quotients is zero \cite{ritter}. 

So from \eqref{eq.FukayaTeleman} we expect the spectral component $\cW_{S^1}\!(X)_\lambda$ to be equivalent to $\cW(T^*S^2)$ for one value of $\lambda$, either $\lambda=1$ or $\lambda=-1$, and zero for all other values of $\lambda$. To nail down the sign one can compute that the Lagrangian $S^3$ in $X$, with its unique spin structure, defines an object of $\cW_{S^1}\!(X)_1$ (see Remark \ref{rem.spinonspheres} below).  So we expect that:
$$\cW_{S^1}\!(X)_1 \cong \cW(T^*S^2) \aand \quad\quad \cW_{S^1}\!(X)_\lambda \cong 0 \;\;\mbox{for}\; \lambda\neq 1$$
 But we can upgrade this, because $X$ is actually preserved by a rank 2 torus $T$, acting on $(x,y)$ and $(z,w)$ separately. This induces a residual $S^1$ action on the quotient $T^*S^2$ which is the same action we considered previously. The obvious prediction is that:
\beq{eq.T*S^3}\cW_T(X)_{(\lambda,\lambda^{-1})} \cong \cW_{S^1}\!(T^*S^2)_\lambda
\aand \quad\quad \cW_T(X)_{(\lambda, \mu)} \cong 0 \;\;\mbox{for}\; \lambda\mu\neq 1
\eeq
 Projecting $X$ to $\C$ via $xy$ exhibits it as a rank 2 algebraic torus fibration degenerating over the two points $D_0=\{0\}, D_1=\{1\}$, as in Conjecture \ref{conj.conics}. Taking the $S^1$ quotient we recover our previous picture of $T^*S^2$ as a rank 1 torus fibration degenerating over the divisor $D=D_0\cup D_1$. 

Now we can construct the toric mirror, and since we have a rank 2 torus action here the mirror should be defined over $\Spec \C[s_1^{\pm}, s_2^{\pm}]$. Following the prescriptions of the previous section the result is the Landau-Ginzburg model $\C^3_{\my, \mw, \mp} \times (\C^*)^2_{s_1, s_2}$ with:
\beq{TS2mirror}\mW = \my\mw\mp + (1+s_1)\my + (1+s_2) \mw \eeq
To get the mirror to $T^*S^2$ we restrict to a divisor $s_1s_2=1$, which recovers \eqref{nodesmoothing2} after an obvious co-ordinate change. It's also easy to check that setting $s_1s_2=\lambda$ for any value $\lambda\neq 1$ results in a superpotential with no critical points. So our prediction \eqref{eq.T*S^3} is consistent with the mirror.

 If we set $s_1=-1$ and $s_2=-1$ we get the mirror to the pair-of-pants again. But with this approach we get a $\C[t_1, t_2]$ structure on $MF(\my\mw\mp)$, the two $HH^2$ classes being $t_1=\my$ and $t_2=\mw$. They are mirror to capping off two punctures separately. 

Also note that $\Crit(\mW)$ is contained in $\{s_1=s_2=-1\}$, which is mirror to the fact that $\cW_T(X)_{(-1,-1)}$ is the only non-zero spectral component.
\end{egu}

The previous example has an immediate higher-dimensional generalization.

\begin{eg} We can get the symplectic manifold $T^*\C\P^n$ as a hyperk\:ahler quotient of $\C^{2n +2}$. We take the hypersurface
$$X  = (x_0y_0 + ... +  x_n y_n  - 1) \;\; \cong\; T^*S^{2n+1} $$
and then take the symplectic quotient of $X$ by the $S^1$ action having weights $(1,-1,,..., 1,-1)$. At $\mma=0$ we get $T^*\C\P^n$ and at other values of $\mma$ we get something diffeomorphic but non-exact. So we expect to see a single non-zero spectral component, and a computation with spin structures (Remark \ref{rem.spinonspheres}) tells us to expect it at $\lambda=(-1)^{n+1}$. All other spectral components should be zero.

Moreover, since this $S^1$ is just a 1-parameter subgroup of a rank $n+1$ torus $T$ acting on $X$, we expect that:
\beq{eq.T*S^odd} \cW_T(X)_{(\lambda_0,..., \lambda_n)} \cong \begin{cases} \cW_{T/S^1}\!(T^*\C\P^n)_{(\lambda_1,..., \lambda_n)} & \mbox{for}\; \lambda_0\lambda_n = (-1)^{n+1} \\
0 & \mbox{for}\; \lambda_0\lambda_n \neq (-1)^{n+1}
\end{cases}
\eeq

We have a projection
 $$ \pi = (x_1y_1,\, ...,\, x_ny_n): \; X \to \A^n$$
which exhibits $X$ as an algebraic torus fibration as in Conjecture \ref{conj.conics}, degenerating over a `simplex' of hyperplanes:
$$D = \bigcup_{i=1}^n \{z_i=0\} \cup \{z_1+...+z_n = 1\}$$
After quotienting by $S^1$ we get a similar picture for $T^*\C\P^n$ (but now this is a slight generalization of the construction in Conjecture \ref{conj.conics}). The complement $\A^n\setminus D$ of this hyperplane arrangement can be viewed as a higher-dimensional analogue of the pair-of-pants. 

The toric mirror to $X$ is $\C^{n+2}_{\mp,\my_0,..., \my_n}\times (\C^*)^{n+1}_{s_0,..., s_1}$ with 
\beq{mirrortoPn}\mW = \mp\my_0...\my_n + \sum_{i=0}^n (1+s_i)\my_i \eeq
where the $(S^1)^n$ action on $X$ is mirror to the given $\C[s_0^{\pm 1}, ...., s_n^{\pm 1}]$ structure. If we set all $s_i=-1$ this reduces to $\MF(\mp\my_0...\my_n)$ which was proven by \cite{LPsym} to be the mirror to $\A^n\setminus D$. So this agrees with our Conjecture \ref{conj.conics}. 
\pgap

To get the mirror to $T^*\C\P^n$ we should restrict to the divisor $\{s_0...s_n = (-1)^{n+1}\}$, this gives:
$$\mW' =  \mp\my_0...\my_n + \left(1+\frac{(-1)^{n+1}}{s_1...s_n}\right)\my_0 +  \sum_{i=1}^n (1+s_i)\my_i  $$

If we're prepared to forget the information of the torus action then we can simplify this considerably. The critical locus of $\mW'$ is contained in $\{s_i=-1, \forall i\}$, so by completing at that locus and repeated applications of \eqref{h2} we can reduce this to the model:
$$\C^2_{\mp,\my}, \quad \mp\my^{n+1} $$
So the claim is that 
$$\cW(T^*\C\P^n) \cong \MF(\mp\my^{n+1})$$
and it would not be very hard to verify this directly.\footnote{The Lagrangian $\C\P^n$ should map to the structure sheaf along $\mp=0$, and the cotangent fibre to the structure sheaf along $\my=0$.}   But this description obscures the $\C[s_1^{\pm 1},..., s_n^{\pm 1}]$ structure. 

\end{eg}

\begin{remu} \label{rem.spinonspheres} An $S^1$ action on $S^k$ induces a Hamiltonian $S^1$ action on $X=T^*S^k$, and then the invariant Lagrangian sphere in $X$ defines an object of either $\cW_{S^1}\!(X)_1$ or $\cW_{S^1}\!(X)_{-1}$ as explained in Remark \ref{rem.signspin}.  To see if it's $+1$ or $-1$ we need a simple computation with spin structures, generalizing the one we did for $T^*S^2$.

 Alternatively, following Remark \ref{rem.exp(a)}, we can ask whether the Lagrangian correspondence $\Gamma$ gives us an object of $\cW_{S^1}\!(X \times X\sslash S^1)_1$ or $\cW_{S^1}\!(X \times X\sslash S^1)_{-1}$. But this is the same computation. 

The frame bundle of $S^k$ is (up to homotopy) just $SO(k+1)$, and the unique spin structure on $S^k$ is the double cover $\operatorname{Spin}(k+1)\to SO(k+1)$.  Suppose our $S^1$ action on $S^k$ is a subgroup of $SO(k+1)$ which rotates a plane and fixes the orthogonal $k-1$ subspace. The induced action on the frame bundle is the obvious one, and orbits are non-trivial elements of $\pi_1(SO(k+1))$. It follows that, for this $S^1$ action, the Lagrangian sphere gives an object of $\cW_{S^1}\!(X)_{-1}$. 

But now suppose that $k=2n+1$ is odd and we take the free $S^1$ action on $S^k$. Since this rotates $n+1$ orthogonal planes the resulting orbits in the frame bundle are contractible if and only if $n+1$ is even, and it follows that the Lagrangian sphere gives an object of $\cW_{S^1}\!(X)_{(-1)^{n+1}}$. 

\end{remu}

We should be able to handle any hyperk\"ahler quotient of a vector space by a torus using the methods of the previous two examples. 

\begin{egu} Let $Y$ be the cotangent bundle to the blow-up of $\C\P^2$ in a point. We can construct $Y$ as a hyperk\"ahler quotient of $\C^8$ by $(U(1))^2$ acting with weights:
 $$\begin{pmatrix} 1 & 1 & 1 & 0 & -1 & -1 & -1 & 0 \\ 0 & 0 & 1 & 1 & 0 & 0 & -1 & -1 \end{pmatrix} $$
 So it's the symplectic quotient of the 6-dimensional affine variety 
$$X = \{x_1y_1 + x_2y_2 + x_3y_3 = \alpha, \;x_3y_3 + x_4y_4 = \beta\} \subset \C^8 $$
(for generic $\alpha, \beta$) by a rank two torus.  Projecting to $\A^2$ with $(x_1y_1, x_2y_2)$ exhibits $X$ as an algebraic torus fibration,  degenerating over a divisor which is the union of four hyperplanes:
$$D = \big\{z_1z_2(z_1+z_2 - \alpha)(z_1 + z_2 + \beta-\alpha) = 0 \big\}$$ 
The toric mirror to $X$ can be constructed using the heuristics of Section \ref{sec.generalheuristics}, the first step is pass to an LG model on $\C^{10}$. The result is the toric variety
\beq{variety}\operatorname{Tot}\big\{\cO(-1)^{\oplus 2}_{\mp,\mq} \to \P^1_{\um:\mv}\big\}\times \C^3_{\my_1,\my_2,\my_3}\times (\C^*)^4_{s_, s_2, s_3, s_4} \eeq
with the superpotential:
$$\mW = \mq\um\my_1\my_2 + \mp\mv\my_3 + (1+s_1)\my_1 + (1+s_2)\my_2 + (1+s_3)\my_3 + (1+s_4)\mq\mv  $$

Now suppose that $\alpha$ and $\beta$ are both real, so that $D$ is the complexification of a real hyperplane arrangement, and suppose further that $\beta>\alpha>0$. Then we can cut $\A^2\setminus D$  into two open sets, one the complement of a triangle of hyperplanes and the other the complement of a square of hyperplanes, \emph{c.f.}~\cite[Example 2.1.5]{Proudfoot2007}. The first open set is topologically  a `4d pair-of-pants', this is the hyperplane arrangement that arises from $T^*\C\P^2$ as we saw in the previous example. The second open set is the cross-product of two copies of the 2d pair-of-pants, this would arise from $T^*\C\P^1\times T^*\C\P^1$. 

On the B-side, the variety \eqref{variety} can be covered by two open charts. Let's set $s_1=s_2=s_3=s_4=-1$, leaving only the terms $\mq\um\my_1\my_2 + \mp\mv\my_3$ in $\mW$. Then on the chart  $\{\mv\neq 0\}$ we can remove the variables $\mp$ and $\my_3$  by heuristic \eqref{h2} and what remains is exactly the mirror to the 4d pair-of-pants. On the chart $\{\um\neq 0\}$ we get the cross-product of two copies of the mirror to the 2d pair-of-pants. So at this level, the decomposition into two open sets can done on either side of the mirror with the same results. 

However if we repeat the B-side discussion without setting the $s_i$'s to -1 then the second chart is not a cross-product, reflecting the fact that $Y$ is not actually the union of $T^*\C\P^2$ and $T^*\C\P^1\times T^*\C\P^1$. 
\end{egu}

The procedure of `cutting into open subsets and regluing' has been used previously to understand mirror symmetry in various situations. The preceding example suggests that mirror symmetry for hypertoric varieties could also be approached this way, and it would be interesting to develop this into a general story.

\subsection{Beyond torus fibrations}\label{sec.beyondconics}

We conclude with two examples which are not algebaric torus fibrations, but do fit into the more general framework of Section \ref{sec.Hamiltonianreduction}, although the manifolds involved are not Liouville. 

\begin{eg} Let $X=\C\P^1\times\C\P^1$ with the diagonal $S^1$ action:
$$(x:y,\, z:w) \mapsto (\lambda x: y,\, \lambda z: w)$$
The moment polytope is an interval, and for a generic value in the interval the symplectic quotient is $\C\P^1$. But there is one singular value in the interior of the interval. The smooth locus of the singular fibre is $U=S^1\times \C^*$, so $U/S^1=\C^*$. 

If we want to think of this as a torus fibration the best we can do is to take the rational map:
$$\pi= xw\!:\!yz \; : \;\C\P^1\times\C\P^1 \dashrightarrow \C\P^1$$
Away from the two base points this is a rank one algebraic torus fibration with two singular fibres.

The mirror to $X$ is the LG model $(\C^*)^2_{\mx, \mz}$ with $\mW= \mx+ 1/\mx + \mz +1/\mz$, and the mirror to the $S^1$ action is $\sigma = \mx\mz$. Changing co-ordinates gives:
$$(\C^*)^2_{\mx, s}, \; \; \mW = (1+s^{-1})\mx + (1+s)/\mx$$
If we restrict to a generic value of $s$ we get the mirror to $\C\P^1$ (with some particular coefficients). But at $s=-1$ we get $\C^*$ with a zero superpotential, as predicted by Conjecture \ref{conj.singularfibres}. 

This example is obvious a compactification of Example \ref{eg.simplest2}, and this mirror is obviously a deformation of that mirror. A slightly more interesting observation is that we can delete the locus $\{xw=yz\}$ and get an honest rank 1 algebraic torus fibration over $\C$ with two critical points, and this is symplectomorphic to $T^*S^2$. So we expect this example to also be a compactification/deformation of Example \ref{eg.TS2again}. On the mirror this is certainly true: we can deform \eqref{nodesmoothing2} to 
$$ \mW = \my\mp(\um - h)  + (1-s)(\my - \mp) $$
and then for a generic value of $h$ this is equivalent to $(\C^*)^2$ with a superpotential $(1+s)(\my - h/\my)$. 

\end{eg}

\begin{egu} Consider the variety $X=K_{\C\P^1}$, which we can get as a symplectic quotient of $\C^3_{x,y,p}$ by $S^1$ acting with weights $(1,1,-2)$. It is diffeomorphic, but not symplectomorphic, to $T^*S^2$. 

Now consider the $S^1$ action on $X$ induced from the action with weights $(-1,0,1)$ on $\C^3$. This preserves the projection:
$$\pi=xyp: X \to \C$$ 
This is not quite the setup of Conjecture \ref{conj.conics} because a generic fibre of $\pi$ is $\C^*$, but the fibre over $0$ is the toric boundary of $X$ and this has two nodes. For a generic value of the moment map the symplectic quotient is $\C$, but there are two singular values, each singular fibre containing one of the nodes.

The mirror to $X$ is $(\C^*)^2$ with 
$$\mW = \mx + \mp + \frac{\mp^2}{\mx}$$
 and the mirror to the $S^1$ action is $\sigma = \mp/\mx$. Changing co-ordinates gives $\mW = (1+s+s^2)\mx$.

So at a generic value of $s$ we get the mirror to $\C$, and at the two roots of the quadratic $1+s+s^2$ we get $\C^*$. This is as predicted by \eqref{eq.FukayaTeleman} and Conjecture \ref{conj.singularfibres}. If we complete $\MF((\C^*)^2, \mW)$ at either of the two special values of $s$ we get a copy of the relative Fukaya category $\cW(\C, \{0\})$, which fits with Remark \ref{rem.multiplesingularfibres}.
\pgap

This example is really a degeneration of $T^*S^2$, the two critical points of the fibration have collided. The hyperk\"ahler perspective of Example \ref{eg.hypertoric1} makes this explicit -  to get $K_{\C\P^1}$ instead of $T^*S^2$ we set the complex moment map to zero and take the symplectic quotient of the singular hypersurface:
$$Y =  (xy + zw) \subset \C^4$$
On the mirror side the fact that $Y$ is singular is not a problem. Removing the constant term from the equation just causes us to delete the divisor $\mp=0$, so the mirror to $Y$ is  $\C^2_{\my,\mw}\times (\C^*)^3_{\mp, s_1,s_2}$ with the superpotential \eqref{TS2mirror}. To get the mirror to $X$ we set $s_1s_2=1$, and after using heuristic \eqref{h2} to remove two variables we recover the mirror we just found.

\end{egu}

The previous example suggests that, when studying mirror symmetry for hypertoric varieties, there should be no problem handing degenerate hyperplane arrangements where multiple hyperplanes coincide.
\pgap

For our final example we consider a case where the symplectic side is not just a manifold $X$, but rather a Landau-Ginzburg model $(X,W)$. It suggests that we can apply all of the ideas of Section \ref{sec.intro} to this more general situation.

\begin{eg} As in Example \ref{eg.simplest0} we take $X=\C^2$ with the $S^1$ action having weights $(1,-1)$. Now add the $S^1$-invariant superpotential $W=xy$. 

The presence of $W$ induces a partial compactification of the mirror \eqref{eq.mirror1} that we found in Example \ref{eg.simplest2}. It's easy to compute that the result is:
$$\big(\C_{\my} \times \C^*_s, \; (1+s)\my\big) $$
Indeed this has one non-degenerate critical point, so the category of matrix factorizations is equivalent to $D^b(pt)$, which matches the Fukaya-Seidel category of $W$. 

If we set $s$ to any value $\lambda\neq -1$ we get $\C_{\my}$ with superpotential $(1+\lambda)\my$ and the category of matrix factorizations of this is zero, but for $\lambda=-1$ we get $D^b(\C)$. So mirror symmetry suggests that the $S^1$-equivariant Fukaya-Seidel categories of $W$ (assuming such things can be defined) will be zero except the for the component at $\lambda =-1$. 

The symplectic quotient of $X$ at a regular value is $\C$, and $W$ descends to give the co-ordinate function. At the singular value our conjectures tell us to replace $\C$ with $\C^*$. The Fukaya-Seidel category of $(\C, z)$ is zero, but the Fukaya-Seidel category of $(\C^*, z)$ is equivalent to $D^b(\C)$ (its toric mirror). So the pictures from mirror symmetry and from symplectic reduction match each other.
\end{eg}
\vspace{20pt}

\subsection*{Acknowledgements.} We have benefited from conversations, suggestions and correspondence from Denis Auroux, Kenji Fukaya, Dan Pomerleano, Ivan Smith, and Jack Smith. YL is partially supported by the Royal Society URF\textbackslash R\textbackslash180024 and EPSRC grant EP/W015889/1.

\bibliographystyle{halphanum}

\end{document}